\documentclass[%
  twocolumn, a4paper
]{mpi2015-cscpreprint}

\usepackage[american]{babel}

\usepackage{graphicx}

\usepackage{amssymb}
\usepackage{amsthm}

\newtheorem{theorem}{Theorem}
\newtheorem{lemma}{Lemma}

\newtheorem{proposition}{Proposition}

\newtheorem{remark}{Remark}

\newcommand{\blem}{\begin{lemma}}
	\newcommand{\elem}{\end{lemma}}
\newcommand{\bea}{\begin{eqnarray*}}
	\newcommand{\eea}{\end{eqnarray*}}

\usepackage{threeparttable} 
\usepackage{amsmath}
\usepackage{multirow}
\usepackage{float}
\usepackage{subfig}
\usepackage[normalem]{ulem}

\usepackage{tikz}
\usepackage{pgfplots}
\pgfplotsset{compat=newest}
\usepackage{xcolor}
\usepackage{algorithm}
\usepackage{algorithmic}
\usepackage{enumerate}
\usepackage{caption}

\begin{document}
  

\title{\emph{Inf-Sup}-Constant-Free State Error Estimator for Model Order Reduction of Parametric Systems in Electromagnetics}

\author[$\ast$]{Sridhar Chellappa}
\affil[$\ast$]{Department of Computational Methods in Systems and Control Theory, Max Planck Institute for Dynamics of Complex Technical Systems, 39106 Magdeburg, Germany.\authorcr
	\email{chellappa@mpi-magdeburg.mpg.de}, \orcid{0000-0002-7288-3880} \authorcr
	\email{feng@mpi-magdeburg.mpg.de}, \orcid{0000-0002-1885-3269} \authorcr
	\email{benner@mpi-magdeburg.mpg.de}, \orcid{0000-0003-3362-4103}}

\author[$\ast$]{Lihong Feng}

\author[$\dagger$]{Valent\'in de la Rubia}
\affil[$\dagger$]{Departamento de Matem\'{a}tica Aplicada a las {TIC}, ETSI de Telecomunicaci\'{o}n, Universidad Polit\'{e}cnica de Madrid, 28040 Madrid, Spain.\authorcr
	\email{valentin.delarubia@upm.es}, \orcid{0000-0002-2894-6813}}

\author[$\ast$]{Peter Benner}

\shorttitle{\emph{Inf-Sup}-Constant-Free State Error Estimator}
\shortauthor{S. Chellappa, L. Feng, V. de la Rubia, P. Benner}
\shortdate{}

\keywords{	Computational techniques, error estimation, finite element method, Galerkin method, microwave circuits, numerical analysis, reduced basis method, reduced order modelling.}

\msc{65M15, 65M60, 65R10, 78A50}
  
\abstract{%
A reliable model order reduction process for parametric analysis in electromagnetics is detailed. Special emphasis is placed on certifying the accuracy of the reduced-order model. For this purpose, a sharp state error estimator is proposed.
Standard \textit{a posteriori} state error estimation for model order reduction relies on the \emph{inf-sup} constant. For parametric systems, the \emph{inf-sup} constant is parameter-dependent. The \textit{a posteriori} error estimation for systems with very small or vanishing \emph{inf-sup} constant poses a challenge, since it is inversely proportional to the \emph{inf-sup} constant, resulting in overly pessimistic error estimation especially at and around resonance frequencies. Such systems appear in electromagnetics since the \emph{inf-sup} constant values are close to zero at points close to resonant frequencies, where they eventually vanish. We propose a novel \textit{a posteriori} state error estimator which avoids the calculation of the \emph{inf-sup} constant. The proposed state error estimator is compared with the standard error estimator and a recently proposed one in the literature. It is shown that our proposed error estimator outperforms both existing estimators. Numerical experiments are performed on real-life microwave devices such as narrowband and wideband antennas, two types of dielectric resonator filters as well as a dual-mode waveguide filter. These examples show the capabilities and efficiency of the proposed methodology.
}

 \novelty{
We propose efficient, sharp a posteriori state error estimation for reduced-order models of general linear parametric systems. Robustness of the estimator is demonstrated via its applications in computational electromagnetics. It avoids computing a quantity that normally tends to blow up existing error estimators at or near resonance frequencies. Our new estimator is integrated within an adaptive greedy algorithm that iteratively
builds the reduced-order model. }

\maketitle
  
\section{Introduction}%
\label{sec:intro}

The rapid rise in demand for novel communication devices with an emphasis on optimal performance has posed a challenge to the standard design and prototyping workflow to manufacture such devices. As a result, a great computational effort is carried out to get physical insight from time-consuming electromagnetic simulations by means of parametric studies. The ultimate goal of these parametric analyses is to develop robust and effective electrical designs, which are of paramount importance for industry. Different efforts in the computational electromagnetic (CEM) community have been carried out to speed this costly process up and many of them follow the model order reduction (MOR) trend \cite{Edlinger2014APosteriori,nicolini2019model,Silveira2014Reduced-OrderModels,Rewienski2016greedy,codecasa2019exploiting,ChewTAP2014GeneralizedModal,xue2020rapid,mrozowski2020,morCheHM11,morHesB13}.

MOR has demonstrated its robustness in reducing the complexity of parametric systems~\cite{morPatR07, morBenGW15}. While many works proposing new MOR methods are available, not all of them guarantee the accuracy of the proposed methods by providing an error analysis of the methods, e.g.,  \cite{nicolini2019model,Silveira2014Reduced-OrderModels}. It is quite usual that the proposition of corresponding error estimation lags behind new algorithms. So far, the standard \textit{a posteriori} state error estimator is widely used if no better choice is available, namely, residual norm divided by the \emph{inf-sup} constant~\cite{morHesB13,morPatR07,hess2015estimating, Edlinger2015CertifiedDualCorrected}. As already stated in~\cite{morFenB19}, keeping the~\emph{inf-sup} constant in the denominator of the error estimator causes risk for many problems with small~\emph{inf-sup} constants. The system of time-harmonic Maxwell's equations is one of them since it has~\emph{inf-sup} constants close to zero near to resonant frequencies. This leads to an error estimator having very large magnitude at those frequencies, even though the true error is already very small~\cite{Pra19,GarRM17,RozV07}. 

Some \textit{a posteriori} error estimators independent of the~\emph{inf-sup} constant are recently proposed in~\cite{morFenB19,morSemZP18}. There, instead of computing the~\emph{inf-sup} constant, additional dual
or residual systems are solved to obtain the error estimators. A state error estimator is proposed in~\cite{morSemZP18}, and an output error estimator is considered in~\cite{morFenB19}. The output error estimators are of interest when the accuracy of the outputs needs to be guaranteed by the reduced-order model (ROM). However, for some cases, e.g., for fast frequency sweeps in electromagnetics as well as in acoustics~\cite{Pra19,GarRM17}, the state vector should be accurately approximated by the ROM. The error estimator in~\cite{morFenB19} is motivated by the aim to control the output error and can hardly be applied to estimate the error of the whole state vector. Therefore, state error estimation is important for generating accurate ROMs in such situations. 

As the central contribution of this work, we introduce a new \textit{a posteriori} error estimator for the state error introduced by the ROM approximation. The proposed estimator avoids computing the~\emph{inf-sup} constant. Instead, it exploits a residual system to construct the estimator. The same residual system was also used in~\cite{morSemZP18}. However, our newly proposed error estimator is more accurate and more efficient to compute than the one detailed in~\cite{morSemZP18}, according to our analyses in Section~\ref{subsec:comp}. To avoid the expensive \emph{inf-sup} constant computation, the norm of the residual introduced by the ROM is often heuristically used as an error estimator. On the one hand, from our analysis in Section II-A, it is not theoretically reliable, and may not guarantee the accuracy of the ROM for certain problems. On the other hand, although it gives acceptable results for many problems, the residual estimator often overestimates the true error, often with effectivity factor much larger than $10$. The proposed error estimator not only theoretically guarantees the accuracy of the ROM, but is also much more accurate, with the effectivity factor around $1$. Further, an adaptive algorithm is proposed to iteratively build the ROM by using this new state error estimator.  

The rest of the paper is organized as follows. In Section~\ref{sec:prob}, we introduce the parametric problem in electromagnetics considered in this work, and the standard state error estimation. The new state error estimation is proposed in Section~\ref{sec:new_erest}. The estimator proposed in~\cite{morSemZP18} is reviewed and theoretically compared with our proposed estimator in Section~\ref{sec:rand_est}. 
Numerical simulations in Section~\ref{sec:numer} show the performance of the proposed estimator, as well as 
the results of the standard estimator and the estimator in~\cite{morSemZP18}. Finally, in Section~\ref{sec:conclusions}, we provide conclusions.

\section{Parametric Problem and Standard State Error Estimation}
\label{sec:prob}

The systems we are interested in are steady (time-harmonic), linear, parametric systems in the form of,
\begin{equation}
\label{eq:fom}
A(\mu) x(\mu) = b(\mu).
\end{equation}
Here, $\mu \in \mathbb C^m$ is the vector of parameters, $A(\mu) \in \mathbb C^{n \times n}$ is the system matrix, $b(\mu) \in \mathbb C^{n}$, and $x(\mu) \in \mathbb C^{n} $ is the state solution, $n$ is assumed to be large, say $n \in \mathcal{O}(10^5)$. System (1) is referred to as the full order model (FOM). Such systems arise from numerical discretization of partial differential equations (PDEs) and integral equations (IEs), such as the time-harmonic Maxwell's equations where differential and integral approaches are both possible.

The ROM can be obtained based on Galerkin projection, 
\begin{equation}
\label{eq:rom}
\hat{A}(\mu) z(\mu)= \hat{b}(\mu),
\end{equation}
where the reduced matrices are given by $\hat{A}(\mu) = V^{T} A(\mu) V \in \mathbb C^{r\times r}$, $\hat{b}(\mu) = V^{T} b(\mu) \in \mathbb C^{r}$ and $V \in \mathbb R^{n\times r}$ is the projection matrix spanning the reduced space. $z(\mu) \in \mathbb C^{r}$ stands for the reduced state vector. $r\ll n$ is the order of the ROM.
The approximate state vector obtained from the ROM is $\hat{x}(\mu) = V z(\mu)$, such that $x(\mu) \approx \hat{x}(\mu)$. 
\begin{remark}
	The system setting and the error estimators presented in the following can be readily extended to systems with multiple inputs \cite{morFenB19}, i.e., the right-hand side is a matrix $B(\mu) \in \mathbb{C}^{n \times p}$. Please also refer to the details in Section~\ref{sec:numer}, where a circuit with two ports is used to test the error estimators.
	The true solution in this case is $X(\mu) \in \mathbb{C}^{n \times p}$ and the corresponding approximate solution matrix is $\hat{X}(\mu) = V Z(\mu), Z(\mu) \in \mathbb{C}^{r \times p}$. For simplicity of explanation in the following sections, we consider a single-input system.
\end{remark}

Next, we introduce the standard state error estimation and point out the role of the \emph{inf-sup} constant. 
\subsection{Standard State Error Estimation}
\label{subsec:standard_errest}
The residual obtained by substituting the approximate solution into the FOM is given by
\begin{equation}
\label{eq:residual}
r(\mu) = b(\mu) - A(\mu) \hat{x}(\mu).
\end{equation}
Then from (\ref{eq:fom}), we obtain,
\begin{equation}
\label{eq:resi2}
r(\mu)= A(\mu) x(\mu) - A(\mu) \hat{x}(\mu),
\end{equation}
so that
\begin{equation}
\label{eq:err-resi}
\| e(\mu) \| :=\, \|x(\mu) -\hat{x}(\mu)\|= \|A(\mu)^{-1} r(\mu)\|.
\end{equation}
The standard \emph{a posteriori} state error estimation is obtained by invoking the sub-multiplicativity property of the operator norm in~(\ref{eq:err-resi}). 
We state the following proposition:
\begin{proposition}
	If $A(\mu)$ is nonsingular, then the norm of the error $e(\mu)$ in (5) is bounded above and below by the norm of the residual $\| r(\mu) \|_{2}$ as
	\begin{equation}
	\frac{1}{\sigma_{\text{max}}} \| r(\mu) \|_{2} \leq \| e(\mu) \|_{2} \leq \frac{1}{\sigma_{\text{min}}} \| r(\mu) \|_{2}
	\label{eq:std_err_bnd}
	\end{equation}
	with $\sigma_{\text{max}}$ and $\sigma_{\text{min}}$ being the maximal and minimal singular values of the matrix $A(\mu)$, respectively.
	\begin{proof}
		From (4), we can show the following upper bound on the residual:
		\begin{align*}
		\| r(\mu) \|_{2} &= \| A(\mu) e(\mu) \|_{2},\\
		& \leq \| A(\mu) \|_{2} \| e(\mu) \|_{2} = \sigma_{\text{max}} \|e(\mu) \|_{2}.
		\end{align*}
		Next, starting once again from (4), we derive an upper bound for the error:
		\begin{align*}
		r(\mu) &= A(\mu) e(\mu),\\
		\implies e(\mu) &= A(\mu)^{-1} r(\mu).
		\end{align*}
		Hence,
		\begin{align*}
		\| e(\mu) \|_{2} 	 &= \| A(\mu)^{-1} r(\mu) \|_{2}\\
		& \leq \| A(\mu)^{-1}\|_{2} \| r(\mu) \|_{2} = \frac{1}{\sigma_{\text{min}}} \| r(\mu) \|_{2}.
		\end{align*}
		Combining the above two bounds, the proposition is shown to be true.
	\end{proof}
\end{proposition}
\vspace{1em}
%

Through this inequality, we define the upper bound for the state approximation error
\begin{equation}
\| e(\mu) \| = \| x(\mu) - \hat{x}(\mu) \| \leq \frac{1}{\sigma_{\text{min}}} \| r(\mu) \| =: \delta(\mu).
\label{eq:err_std}
\end{equation}
Further, in this discretized setting, $\sigma_{\text{min}}$ plays the role of the \emph{inf-sup} constant.
However, the standard error estimator $\delta(\mu)$ approaches infinity when the matrix $A(\mu)$ is close to singular at some values of $\mu$, resulting in a rather poor estimation and a rough bound. This is true for many problems~\cite{morFenB19}. Furthermore, the above error estimation leads to unacceptable overestimation of the state error even for well-conditioned problems~\cite{morSchWH18}. Given no better choices, the residual norm $\| r(\mu)\|$ is often used as a heuristic error estimator~\cite{Rub14, Rub18, morRewLM15, de2009reliable, Vouvakis2011FastFrequency, Edlinger2015ANewMethod, Edlinger2017finite, fotyga2018reliable, MonjedelaRubia2020EFIE}. However, it is clear from (\ref{eq:std_err_bnd}) that using $\| r(\mu) \|$ as the error estimator is not reliable, as the scaling factor before it is completely ignored. In particular, when $\sigma_{\text{max}} \ll 1$, we have  $\| e(\mu) \| \gg \| r(\mu) \|$, meaning $\| r(\mu) \|$ may \emph{underestimate} the true error much. In the next section, we propose a new error estimator that avoids computing the \emph{inf-sup} constant while being theoretically more reliable than the heuristic estimator $\| r(\mu) \|$.

\section{Proposed State Error Estimation}
\label{sec:new_erest}
From~(\ref{eq:resi2}), we know that
\begin{equation*}
A(\mu)(x(\mu)- \hat{x}(\mu))=r(\mu).
\end{equation*}
To compute the error $e(\mu)$ for each value of $\mu$, we need to solve the residual system corresponding to each $\mu$:
\begin{equation}
\label{eq:sys_resi}
A(\mu)e(\mu)=r(\mu).
\end{equation}
Note that the residual system~(\ref{eq:sys_resi}) is of the original large dimension $n$. It is not practical to 
solve many large systems at many values of $\mu$ in order to know the error distribution in the whole parameter domain. For fast error estimation, we first construct the ROM of the residual system using a projection matrix $V_{\text{e}} \in \mathbb{R}^{n \times \ell}$ that spans the error subspace. This is given as,  
\begin{equation}
\label{eq:rom_sys_resi}
\tilde {A}(\mu) z_e(\mu) = \tilde{r}(\mu),
\end{equation}
where $\tilde {A}(\mu) = V_{\text{e}}^{T} A(\mu)  V_{\text{e}}$, $\tilde{r}(\mu) = V_{\text{e}}^{T} r(\mu)$, and $\tilde e(\mu):=V_{\text{e}}z_e(\mu)$ approximates $e(\mu)$.
Then, our proposed state error estimation is given by $\|\tilde e(\mu)\|$, i.e. 
\begin{equation}
\label{eq:err_prop}
\|x(\mu)- \hat{x}(\mu)\|=\|e(\mu)\|\approx \|\tilde e(\mu)\|.
\end{equation}
We obtain the following analysis for the rigorousness and tightness of the state error estimator $\|\tilde e(\mu)\|$.
\begin{theorem}
	\label{thm:errestm_rigor}
	The state error $\|e(\mu)\|$ can be bounded by $\|\tilde e(\mu)\|$ as follows:
	\begin{equation}
	\label{eq:err_bound}
	\|\tilde e(\mu)\|- \gamma(\mu)\leq \|e(\mu)\| \leq \|\tilde e(\mu)\| + \gamma(\mu).
	\end{equation}
	Here $\gamma(\mu)=\|e(\mu)-\tilde e(\mu)\|\geq 0$ is independent of the \emph{inf-sup} constant and it is a small value if $\tilde e(\mu)$ well approximates $e(\mu)$, which is achievable by accurate MOR of the residual system~(\ref{eq:sys_resi}).
\end{theorem}
The theorem can be easily proved by applying the triangular inequality to $\|e(\mu)-\tilde e(\mu)\|$, and is not detailed here.
\begin{remark}
	Theorem \ref{thm:errestm_rigor} is also valid for a complex matrix $V \in \mathbb{C}^{n \times r}$. Considering that a real matrix $V$ is used in our numerical tests (see Remark 5), we keep $V$ real all through the paper to avoid confusion.
\end{remark}
\subsection{Computing the State Error Estimator}
\label{subsec:comp_errest}
In order to compute the state error estimator $\|\tilde e(\mu)\|$, we need to construct the ROM~(\ref{eq:rom_sys_resi}) for the residual system~(\ref{eq:sys_resi}). Consequently, the projection matrix $V_{\text{e}}$ has to be computed. In a similar manner as in~\cite{morFenB19}, we look at the residual system~(\ref{eq:sys_resi}) in order to identify the subspace for the trajectory of the state error vector $e(\mu)$:
\begin{equation}
\label{eq:trj_e}
\begin{array}{rcl}
e(\mu)&=&A(\mu)^{-1}r(\mu)\\
&=&A(\mu)^{-1}(b(\mu) - A(\mu) \hat{x}(\mu))\\
&=&A(\mu)^{-1}b(\mu) - \hat{x}(\mu)\\
&=&A(\mu)^{-1}b(\mu) - Vz(\mu).
\end{array}
\end{equation}
Note that $A(\mu)^{-1}b(\mu)$ is nothing but the state solution $x(\mu)$. Suppose there exists a projection matrix $V_{r} \in \mathbb{R}^{N \times r}$ such that its columns span a subspace $\mathcal{V}_{r}$ in which $x(\mu)$ can be well approximated, then $x(\mu)$ can be approximately represented by the columns of $V_r$, i.e. $x(\mu)\approx V_r z_r(\mu)$ with $z_{r}(\mu) := V_{r}^{T} x(\mu)$. Thus, 
\begin{equation}
\label{eq:trj_e2}
\begin{array}{rcl}
e(\mu)&=&A(\mu)^{-1}b(\mu) - Vz(\mu)\\
&\approx&V_r z_r(\mu)-Vz(\mu).
\end{array}
\end{equation}
It is clear from~(\ref{eq:trj_e2}) that $e(\mu)$ can be approximated by  the linear combination of the columns in $V_r$ and $V$. Therefore, we can construct $V_e$ as 
\begin{equation}
\label{eq:Ve}
V_e=\textsf{orth}([V_r,\, V]),
\end{equation}
where we use MATLAB\textsuperscript \textregistered \hspace{0.01pt} notation to denote the orthogonalization of the column space of the input argument of the \textsf{orth} function.

Since the state error vector $e(\mu)\neq 0$, $\textrm{range}(V_r)$ should be different from $\textrm{range}(V)$  used for constructing the original system ROM~(\ref{eq:rom}). Notice that, if $V_{r} = V$, this implies $V_{e} = V$. The following theorem holds:
\begin{theorem}
	If $V_{e} = V$, then $\| \tilde{e}(\mu) \| = 0$.
	\begin{proof}
		Suppose that $V_{e} = V$. The ROM (9) reads
		\[
		\hat{A}(\mu) z_{e}(\mu) = \tilde{r}(\mu),
		\]
		where we have used the fact that $\tilde{A}(\mu) = \hat{A}(\mu)$ for $V_{e} = V$. Substituting (3) into $\tilde{r}(\mu) = V^{T} r(\mu)$, we get
		\begin{align*}
		\tilde{e}(\mu) &= V z_{e}(\mu),\\
		&= V \big( \hat{A}(\mu)^{-1} (V^{T} r(\mu))  \big),\\
		&= V \big( \hat{A}(\mu)^{-1} V^{T} (\,b(\mu) - A(\mu) \hat{x}(\mu)\,)\big)\,\, (\textnormal{using}\, (3)),\\
		&= V \big(  \hat{A}(\mu)^{-1} (\,\hat{b}(\mu) - \hat{A}(\mu) z(\mu)\,)\big),\\
		&= V \big( \hat{A}(\mu)^{-1} (0) \big)\,\, (\textnormal{using}\, (2)),\\
		&= 0.
		\end{align*}
	\end{proof}
\end{theorem}
The computation of $V_{r}$ (avoiding $V_{r} = V$) will be detailed in Algorithm 1.

In the next subsection, we propose an algorithm for constructing the ROM~(\ref{eq:rom}) of the original FOM~(\ref{eq:fom}) by using the proposed state error estimator $\|\tilde e(\mu)\|$. 

\begin{remark}
	The expression for the error $e(\mu)$ in (13) is seen to be similar to a heuristic approach to check the accuracy of a ROM by using two separate reduced-order models of varying orders and determining the difference between their approximation of the state vector. However, we emphasize that our approach has its differences. Our approach uses $\tilde{e}(\mu)$ defined in (10), which has a guaranteed accuracy as per Theorem 1. Whereas, using the norm of $V_{r} z_{r}(\mu) - V z(\mu)$ in (13) as an error estimator is not guaranteed to be accurate. Moreover, the heuristic approach does not provide any guidance on which criteria to use to compute $V_{r}$ and $V$ so that $\|V_{r} z_{r}(\mu) - V z(\mu)\|$ is an accurate estimation.  Therefore, this approach can only be heuristically implemented.
\end{remark}

\subsection{Constructing the ROM}
The algorithm for constructing the ROM~(\ref{eq:rom}) is detailed in Algorithm~\ref{alg:rom_errest}, where the proposed state error estimator is used to select samples of the parameter $\mu$ for computing the projection matrix $V$.  
\begin{algorithm}[t!]
	\caption{Constructing the ROM~(\ref{eq:rom}) using the state error estimator $\|\tilde e(\mu)\|$}
	\label{alg:rom_errest}
	\begin{algorithmic}[1]
		\REQUIRE System matrix and right-hand side vector $A(\mu), b(\mu)$, training set $\Xi$ including a certain number of samples of $\mu$, tolerance $\texttt{tol}$ for the acceptable state error.
		\ENSURE ROM~(\ref{eq:rom}).
		\STATE Initialize $V = [\,]$, $V_{e} = [\,]$, $\epsilon = 1 + \texttt{tol}$.  Two different samples $\mu^*$ and $\mu_e^*$ randomly taken from $\Xi$.
		\WHILE{$\epsilon > \texttt{tol}$}
		\STATE Compute $V(\mu^*)$ using a favorite MOR method and update $V$:
		$V = \textsf{orth}([V, V(\mu^*)])$.
		\STATE Compute $V_r(\mu_e^*)$ using a favorite MOR method and update $V_r$: 		$V_{r} = \textsf{orth}([V_r, V_r(\mu_e^*)])$.
		\STATE If a real $V$ is preferable, then\\ $V := \textsf{orth}\big(\big[\text{real}(V)\,,\,\text{imag}(V)\big]\big)$; also, if a real $V_{r}$ is preferred, then $V_{r} := \textsf{orth}\big(\big[\text{real}(V_{r})\,,\,\text{imag}(V_{r})\big]\big)$.
		\STATE Form $V_e$: $V_{e}=\textsf{orth}([V_r,\, V])$.
		\STATE Compute the estimated state error vector $\tilde e(\mu)$ using $V$ and $V_e$.
		\STATE Choose the next sample $\mu^*$ from $\Xi$ as $$\mu^*=\arg\max\limits_{\mu\in \Xi} \|\tilde e(\mu)\|.$$
		\STATE Choose the next sample $\mu_e^*$ from $\Xi$ following $$\mu_e^*=\arg\max\limits_{\mu\in \Xi} \|r_e(\mu)\|,$$ where $r_e(\mu)=r(\mu)-A(\mu)\tilde e(\mu)$.
		\STATE $\epsilon=\|\tilde e(\mu^*)\|.$
		\ENDWHILE
		\STATE Use $V$ to construct the ROM: $\hat A(\mu)=V^TA(\mu)V$, $\hat b(\mu)=V^T b(\mu)$.
	\end{algorithmic}
\end{algorithm}
We make some remarks to highlight various aspects of Algorithm~\ref{alg:rom_errest}.
\begin{itemize}
	\item $\textsf{orth}(\cdot)$ in Algorithm 1 means to orthonormalize the columns of the matrices inside the parentheses to get a single orthonormal matrix. This can be done by, e.g., modified Gram-Schmidt process or QR decomposition.
	\item The algorithm is automatic. The user needs only to provide a training set $\Xi$. Adaptive sampling approaches for iteratively updating the training set exist~\cite{morEftPR10,morHesZ16,morJiaCN17,morCheFB19} and can be combined with Algorithm~\ref{alg:rom_errest}. Since it is not the focus of the paper, we will present the corresponding algorithm elsewhere.
	\item In Steps 3-4, parameter dependent matrices $V(\mu^*)$ and $V_r(\mu_e^*)$ 
	can be computed using a favorite MOR method. When using the reduced basis method, we have 
	$$V(\mu^*)=A(\mu^*)^{-1}b(\mu^*), \ V_r(\mu_e^*)=A(\mu_e^*)^{-1}b(\mu_e^*).$$ In~\cite{morFenB19}, a multi-moment-matching method~\cite{morFenB14} is used to compute $V(\mu^*)$, which can also be used to compute $V_r(\mu_e^*)$. We do not repeat the details in this work.
	\item It is important that the parameter samples $\mu^*$ and $\mu_e^*$ must be different, otherwise, 
	this may lead to $V=V_r$, which is in contradiction with the analysis in Subsection~\ref{subsec:comp_errest} (see~(\ref{eq:trj_e2})). Therefore, we choose two different initial samples, and select two different sequential ones in Step 8 and Step 9, respectively. In Step 9, we simply use $\|r_e(\mu)\|$, the residual norm introduced by $\tilde e(\mu)$, the approximate solution to the residual system~(\ref{eq:sys_resi}), as the indicator for selecting $\mu^*_e$ in a greedy strategy. We have $r_e(\mu)=r(\mu)-A(\mu)\tilde e(\mu)$.
	\item $\epsilon$ is taken as the maximal value of the state error estimator, i.e. $\|\tilde e(\mu^*)\|$.
	\item Algorithm 1 involves increased computational effort (one additional FOM solution at every greedy iteration), when compared to a greedy algorithm using the residual norm $\| r(\mu) \|$ as a heuristic error estimator. However, this is a small price to pay considering the by far increased reliability of the proposed state error estimator $\| \tilde{e}(\mu) \|$. A residual norm can easily over- or underestimate the true error. On the one hand, overestimation of the true error often leads to a ROM whose reduced dimension $r$ is unnecessarily large for the desired tolerance; this is the case when the \emph{inf-sup} constant $\sigma_{\text{min}}$ is large. On the other hand, underestimation of the true error can falsely trigger an early termination of the greedy algorithm, thus resulting in a poor ROM. This is true for systems for which $\sigma_{\text{max}}$ is small, e.g., $\sigma_{\text{max}} \leq 0.1$; see (\ref{eq:std_err_bnd}). We also emphasize that the additional computational cost for the proposed state error estimator is restricted to the \emph{offline stage} of MOR and that there is no reduction of the computational speedup at the online stage.
\end{itemize}
\section{Review of the Randomized State Error Estimator}
\label{sec:rand_est}
In this section, we briefly summarize an existing approach in literature~\cite{morSemZP18}, which also proposes an \emph{inf-sup} constant free \emph{a posteriori} error estimator. We also compare and contrast the key differences between the error estimation methodology from~\cite{morSemZP18} and our proposed method.

A randomized state error estimator is proposed in~\cite{morSemZP18} for Galerkin projection based MOR. It is stated that the $\|\cdot\|_\Sigma$ norm of the state error satisfies
\begin{equation} 
\label{eq:err_rand}
\|e\|_\Sigma^2=e^T\Sigma e=e^T\mathbb E(zz^T)e=\mathbb E((z^Te)^2),
\end{equation}
where $z \in \mathbb R^n$ is a zero mean Gaussian random vector with covariance matrix $\Sigma \in \mathbb R^{n \times n}$ and $\mathbb{E}(\cdot)$ refers to the expected value or the mean value. 

An approximation to $\mathbb E((z^Te)^2)$ is firstly proposed,
\begin{equation} 
\label{eq:err_approx}
\mathbb E((z^Te)^2)\approx \big(\frac{1}{K}\sum\limits_{i=1}^{K}(z_i^Te)^2 \big)^{1/2},
\end{equation}
where $z_i \in \mathbb R^n$ are K samples of $z$.
After simple derivations, it is proved that~\cite{morSemZP18} 
\begin{equation}
\label{eq:err_approx_equi}
\big(\frac{1}{K}\sum\limits_{i=1}^{K}(z_i^Te)^2 \big)^{1/2}=\big(\frac{1}{K}\sum\limits_{i=1}^{K}( \xi_i(\mu)^Tr(\mu))^2 \big)^{1/2},
\end{equation}
where $r(\mu)$ is the residual defined in~(\ref{eq:residual}) and the parametric functions $\xi_i(\mu), i=1,\ldots, K$, satisfy the following $K$ random dual systems, respectively,
\begin{equation}
\label{eq:dual_rand}
A(\mu)^T\xi_i(\mu)=z_i,\qquad 1\leq i\leq K.
\end{equation}
To define the error estimator, reduced systems for those $K$ random dual systems are first constructed as,
\begin{equation*}
V_{rd}^TA(\mu)^TV_{rd}\hat \xi_i(\mu)=V_{rd}^Tz_i,\qquad 1\leq i\leq K.
\end{equation*}
The subscript $rd$ is short for \emph{random dual}. Then $\xi_i(\mu)$ can be approximated by $V_{rd}\hat \xi_i(\mu)$ as $\xi_i(\mu)\approx \tilde \xi_i(\mu):=V_{rd}\hat \xi_i(\mu)$. Finally, the state error estimator is defined as
\begin{equation}
\label{eq:estrand}
\tilde \Delta(\mu)=\big(\frac{1}{K}\sum\limits_{i=1}^{K}(\tilde \xi_i(\mu)^Tr(\mu))^2 \big)^{1/2}.
\end{equation}
It is proved in~\cite{morSemZP18} that if $V_{rd}=V_e$, then $\tilde \Delta(\mu)$  in~(\ref{eq:estrand}) can be written as
\begin{equation}
\label{eq:estrand_var}
\tilde \Delta(\mu)=\big(\frac{1}{K}\sum\limits_{i=1}^{K}(z_i^T \tilde e(\mu))^2 \big)^{1/2}.
\end{equation}
Note that $V_e$ and $\tilde e(\mu)$ are 
defined in~(\ref{eq:rom_sys_resi}). This means, with the assumption
$V_{rd}=V_e$, $\tilde \Delta(\mu)$ is an average of K inner products of $z_i$ and $\tilde e(\mu)$, $i=1,\ldots, K$.

\subsection{Comparison between $\tilde \Delta(\mu)$ and $\|\tilde e(\mu)\|$} 
\label{subsec:comp}
\begin{itemize}
	
	\item It is clear that $\tilde \Delta(\mu)$ is an approximation to\\ $\big(\frac{1}{K}\sum\limits_{i=1}^{K}( \xi_i(\mu)^Tr(\mu))^2 \big)^{1/2}$ in~(\ref{eq:err_approx_equi}). The approximation is caused by MOR for the $K$ random dual systems in~(\ref{eq:dual_rand}); or alternatively by MOR for the residual system~(\ref{eq:sys_resi}) if $\tilde \Delta(\mu)$ in~(\ref{eq:estrand_var}) is considered. Combining~(\ref{eq:err_rand}), (\ref{eq:err_approx}) and~(\ref{eq:err_approx_equi}), we see that\\ $\big(\frac{1}{K}\sum\limits_{i=1}^{K}( \xi_i(\mu)^Tr(\mu))^2 \big)^{1/2}$ is again an approximation to the state error $\|e\|_\Sigma$. Whereas, our proposed state error estimator $\|\tilde e\|$ is a direct approximation to the state error $\|e\|$, where the approximation is only introduced by MOR for the residual system~(\ref{eq:sys_resi}). If similar accuracy of MOR for both error estimators is obtained, then $\|\tilde e\|$ should be more accurate than $\tilde \Delta(\mu)$.
	
	\item For computing $\tilde \Delta(\mu)$, no efficient method is proposed in~\cite{morSemZP18} to compute the projection matrix $V_e$ for the reduced residual system~(\ref{eq:rom_sys_resi}). Instead, two algorithms are proposed \cite[Algorithms 3.1, 3.2]{morSemZP18} to compute $V_{rd}$ and reduce the $K$ random dual systems (\ref{eq:dual_rand}). In contrast, we have proposed an efficient method of computing $V_e$ in subsection~\ref{subsec:comp_errest}, so that only one system needs to be reduced.
	
	\item If $\|\cdot\|_2$ is used for both error estimators, then $\Sigma$ for $\|e\|_\Sigma$ is the identity matrix, according to~(\ref{eq:err_rand}). 
\end{itemize}

\newlength\fheight
\newlength\fwidth
\setlength\fheight{3cm}
\setlength\fwidth{3cm}
\section{Numerical Results}
\label{sec:numer}
We test the proposed state error estimator $\|\tilde e(\mu)\|$ and the existing ones on six different real-life applications. The first is a dual-mode circular waveguide filter, and two excitation ports are considered. The next three examples are models of narrowband and wideband antennas, where only one excitation port is taken into account. To illustrate the wide applicability of the proposed approach, we choose challenging models for the last two examples. The fifth example is that of a coax-fed dielectric resonator filter, consisting of many resonances over the frequency range of interest, while for the final example, we consider a three-parameter inline dielectric resonator filter in which the dielectric constants in two spatial regions of the filter are considered as additional parameters. After numerical discretization of the time-harmonic Maxwell's equations by means of the finite element method (FEM), the systems can be written in the form of (\mbox{\ref{eq:fom}}). The in-house code for FEM simulations uses a second-order first family of N\'ed\'elec's elements \cite{Ned80, Ing06} on meshes provided by \texttt{Gmsh} \cite{GeuR09}. Fast frequency sweep is considered for the first five examples. As a result, only one parameter models are taken into account with $\mu = s := \jmath 2 \pi f$. Here, $s$ is the complex frequency variable, $\jmath$ is the imaginary unit, and $f$ is the frequency with unit Hz. The system matrix has the affine form given by,
$$A(s) := \mathcal{S} + s \mathcal{U} + s^{2} \mathcal{T},$$
where $\mathcal{S}$ is the stiffness matrix, $\mathcal{T}$ is the mass matrix and $\mathcal{U}$ is the FEM matrix related to the first-order absorbing boundary conditions (ABC) and $\mathcal{S, U, T} \in \mathbb{R}^{n \times n}$. $B(s) := s \mathcal{Q}$ with $\mathcal{Q} \in \mathbb{R}^{n \times p}$, with $p$ being the number of ports and $\mathcal{Q}$ being a matrix related to the excitation currents at the ports. Finally, the state solution $X(s) \in \mathbb{C}^{n \times p}$ stands for the electric field in the analysis domain. It should be pointed out that integral equation methods for electromagnetic scattering can be taken into account as parametric problems in (\ref{eq:fom}) \cite{MonjedelaRubia2020EFIE,wu2019MLACE}. By the same token, other parameters than frequency are also possible \cite{dang2017quasi}. Therefore, for the last example, in addition to the frequency response, we are also interested in the variation of the system response with respect to the two dielectric constants $d_{1}, d_{2}$. We have  $\mu = (s, d_{1}, d_{2})$. The system matrix for the final example has the following affine form:
$$A(s) := \mathcal{S} + s \mathcal{U} + s^{2} \big(\mathcal{T}_{0} + \frac{d_{1}}{d_{\text{ref}}} \mathcal{T}_{1} + \frac{d_{2}}{d_{\text{ref}}}\mathcal{T}_{2}\big)$$
with $\mathcal{T}_{0}, \mathcal{T}_{1}, \mathcal{T}_{2} \in \mathbb{R}^{n \times n}$, and $d_{\text{ref}}$ being the reference dielectric constant of the dielectric material in the analysis domain.

%

In our experiments, we define the true error ($\epsilon_{\text{true}}$) based on the number of input ports. We first consider the error for each column $X_i(\mu)$ of the solution $X(\mu)$ separately and define the \textit{maximal} true error as,
\begin{equation*}
\epsilon_{\text{true}} = \max_{\substack{\\i \in \{1, \ldots, p\} \\ \mu \in \Xi }} \| X_{i}(\mu) - \hat{X}_{i}(\mu) \|,
\end{equation*}	
where $\hat{X}_{i}(\mu)$ is the $i$-th column of the approximate solution $\hat X(\mu)$. Here, $\Xi$ denotes the training set consisting of elements sampled from the parameter space. Note that for each $i$, the system of $X_i(\mu)$ is a single input system, i.e. 
$A(\mu)X_i(\mu)=B_i(\mu)$,
%
where $B_i(\mu)$ is the $i$-th column of $B(\mu)$. 
We use the \textit{maximal error estimator} to estimate $\epsilon_{\text{true}}$: 
\begin{equation*}
\epsilon_{\text{est}} = \max_{\substack{\\i \in \{1, \ldots, p\} \\ \mu \in \Xi }} \Delta_{i}(\mu).
\end{equation*}
$\Delta_{i}(\mu)$ refers to any of the four error estimators used to estimate the error $\| X_{i}(\mu) - \hat{X}_{i}(\mu) \|$: the residual norm $\|r(\mu)\|$ (\ref{eq:residual}), the standard error estimator $\delta(\mu)$ (\ref{eq:err_std}), the randomized error estimator $\tilde{\Delta}(\mu)$ (\ref{eq:estrand_var}), and the proposed state error estimator $\|\tilde{e}(\mu)\|$ (\ref{eq:err_prop}). We use the metric of effectivity to gauge how close the estimated error is to the true error: 
\mbox{
	$	\texttt{eff} := \frac{\text{Error estimator}}{\text{True error}} = \frac{\epsilon_{\text{est}}}{\epsilon_{\text{true}}}$.}
For each of the six examples considered, we evaluate the performance of the four error estimators. In the sequel,
\begin{itemize}
	\setlength\itemsep{0.7em}
	\item Test 1 refers to the greedy algorithm with the standard error estimator $\delta(\mu)$,
	\item Test 2 denotes the greedy algorithm employing the residual norm $\|r(\mu)\|$ as a heuristic error estimator,
	\item Test 3 uses the randomized error estimator $\tilde{\Delta}(\mu)$ from \cite{morSemZP18}, to drive the greedy algorithm and
	\item Test 4 consists of the proposed error estimator $\| \tilde{e}(\mu) \|$ in the greedy algorithm.
\end{itemize}
In all numerical tests, the vector 2-norm $\|  \cdot \|_{2}$ is used for both the true and estimated errors.
\begin{remark}
	Note that the error $\epsilon_{\text{true}}$ and error estimators are defined for the \textit{scaled} solution after scaling the right hand side matrix $B(\mu)$ with a proper scaling constant in order to avoid large norm of the solution $X(\mu)$. This is due to the fact that the entries in the left hand side matrix $A(\mu)$ have much smaller magnitudes than those in $B(\mu)$ because of the large value of $s$ associated with $\mathcal{Q}$. Without scaling, both $X(\mu)$ and $\hat X(\mu)$ have large norms, leading to large absolute errors, which cannot reflect the real accuracy of the ROM. Therefore, we first scale the right hand side and then construct the  ROM. The scaling constant $scale$ is determined by looking at the magnitude of the largest entry in the right-hand side matrix $B(\mu) := s \mathcal{Q}$. The approximate solution $\hat X(\mu)$ can be recovered by scaling back without losing any accuracy. The scaling constant we take here is $\textnormal{scale}=10^{5}$, and the right hand side after scaling is $s \mathcal{Q} / \textnormal{scale}$.
\end{remark}
\begin{remark}
	Algorithm 1 is tailored for the proposed error estimator $\|\tilde e\|$, which can be seen from Steps 7-8 of the algorithm. Similar greedy algorithms can be developed for the other three error estimators. To avoid redundancy, we do not list the corresponding algorithms for those estimators. However, we should point out that although we use the same name: \textit{greedy algorithm} for all the error estimators, they are not the same algorithm, but the \textit{corresponding algorithm} for each of the estimators. The fact that greedy algorithms may be different with different error estimators is mainly due to the fact that different computational steps (e.g., projection matrices) could be involved in computing the error estimators. For example, instead of $V_e$, $V_{rd}$ is needed for computing the randomized error estimator $\tilde \Delta(s)$. In particular, we use the greedy algorithm proposed in~\cite{morSemZP18} for the randomized estimator.
\end{remark}
\begin{remark}
	For the randomized error estimator $\tilde{\Delta}(\mu)$ in Test 3, we make use of Algorithm 3.1 from \cite{morSemZP18} to construct $V_{rd}$, where a separate greedy algorithm is used. For this purpose, $\texttt{tol}_{rd}$ is defined to be the tolerance for this separate greedy procedure to generate $V_{rd}$.
\end{remark}
\begin{remark}
	For the three antenna examples in this work, the state vectors are complex, so that the projection matrix $V$ obtained in Steps 3-4 of Algorithm 1 is complex. However, the system matrices $\mathcal{S},\, \mathcal{U},\, \mathcal{T}$ are all real. It is seen that $\text{colspan}\{V\} \subset \text{colspan}\{\text{Re}(V),\,\text{Im}(V)\}$ over $\mathbb{C}$. Here $\text{colspan}\{\cdot\}$ means the subspace spanned by the columns of a matrix (the matrices). Then, from the proof of Lemma 6.2 in \cite{morFenB14}, it can be easily proved that the ROM constructed using $\widetilde{V} := \text{orth}\big[\text{Re}(V),\,\text{Im}(V)\big]$ and the ROM computed using $V$ satisfy the same moment-matching property in Theorem 6.1 from \cite{morFenB14}. Based on this observation, we redefine a real projection matrix $V$ as $V := \widetilde{V}$ to make the reduced matrices $\hat{\mathcal{S}},\, \hat{\mathcal{U}},\, \hat{\mathcal{T}}$ also real, but still keep the same approximation property as the original complex-valued $V$.
\end{remark}
\vspace{\baselineskip}

Owing to the large size of the models, the simulations for all the examples were performed on a workstation with 3 GHz Intel\textsuperscript \textregistered Xeon E5-2687W v4 processor and 256 GB of RAM, with MATLAB\textsuperscript \textregistered 2019b. It should be pointed out that, in our experiments, the CPU time performance is not optimized, since MATLAB\textsuperscript \textregistered code is used in the numerical computations. However, a fair comparison among the different methodologies is taken into account.

\subsection{Example I : Dual-Mode Circular Waveguide Filter}
\begin{figure}[tbp]
	\centering
	\includegraphics[width=\linewidth]{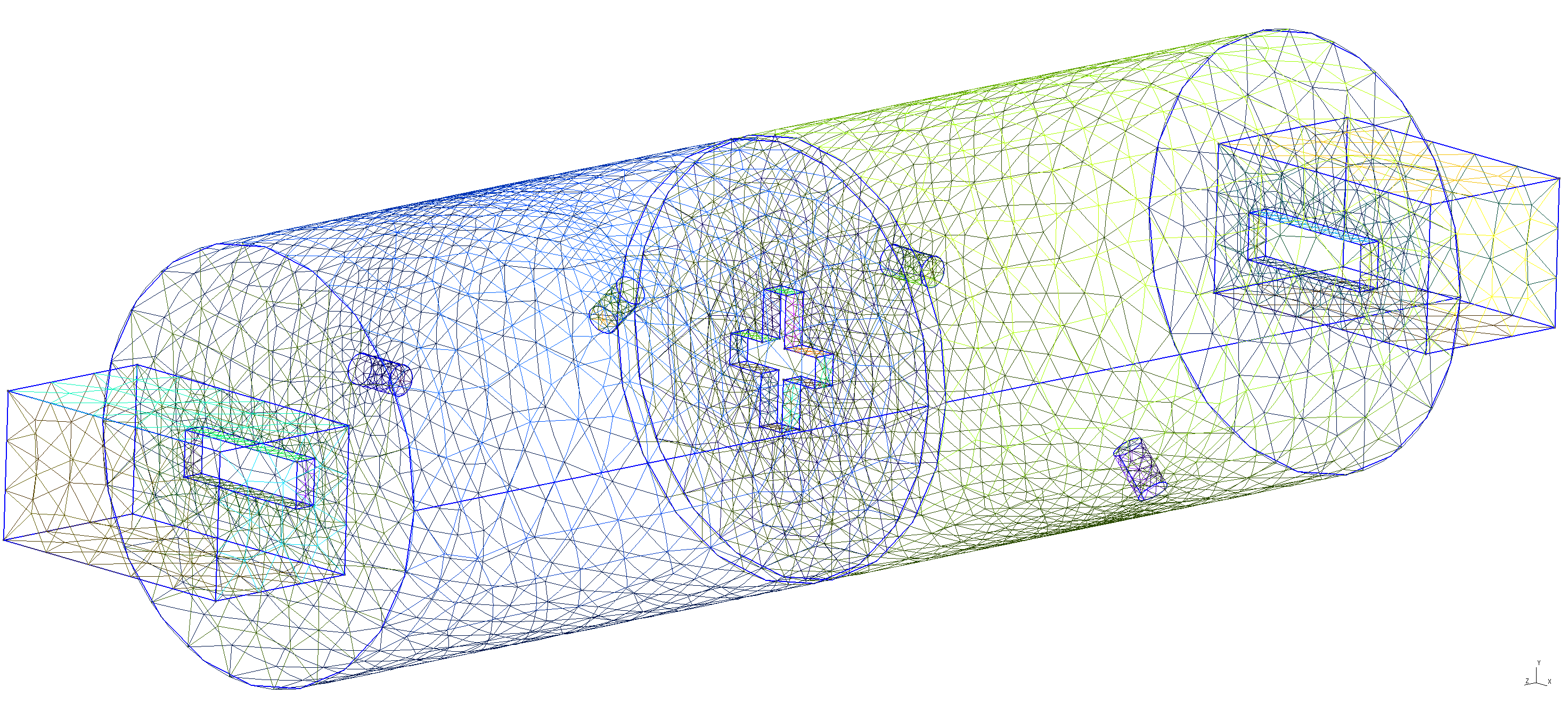}
	\caption{Dual-mode circular waveguide filter. Cavity length 43.87 mm, radius 14 mm, iris thicknesses 1.5 mm, slot lengths 10.05 mm, slot widths 3 mm, arm widths 2 mm, horizontal arm length 7.65 mm, vertical arm	length 8.75 mm, tuning screw depth 3.59 mm and coupling screw depth 3.31 mm.}
	\label{fig:DualModeFilter_Geometry}
\end{figure}

Fig.~\ref{fig:DualModeFilter_Geometry} shows a dual-mode circular waveguide filter as well as its geometry dimensions and mesh for FEM analysis. Dual-mode filters are widely used in satellite communication, due to their power handling capabilities \cite{Rub14, RubZ07}. We consider a system with $p = 2$  input ports. The dimension of the system is $n = 36,426$ and the system matrix is $A(s) = \mathcal{S} + s^{2} \mathcal{T}$. This time, we solve for the E-wall cavity problem and no ABC is needed, as a result, $\mathcal{U} \equiv 0$. The frequency band of interest is [$11.5, \, 12$] GHz. The tolerance for the state error of the ROM is set to be $\texttt{tol}=10^{-6}$. The training set $\Xi$ for Algorithm~\ref{alg:rom_errest} is made up of $101$ uniformly sampled $f\in$ [$11.5, \, 12$] GHz.
\subsubsection{Test 1. Standard error estimator}
For the first test, we use the standard error estimator $\delta(s)$, which involves computing the \textit{inf-sup} constant. Calculating the \textit{inf-sup} constant for all parameters in the training set is expensive. Fig.~\ref{fig:std_conv}  shows that the greedy algorithm requires 10 iterations to converge to the defined tolerance. Fig.~\ref{fig:std_eff}  plots the effectivity changing with iteration number. It is of order $\mathcal{O}(100)$, showing that $\delta(s)$ is not sharp. The algorithm results in a ROM of size 20.\\
\subsubsection{Test 2. Residual norm as a heuristic error estimator}
We use the norm of the residual $\|r(s)\|_2$ as an estimator for the state error. Fig.~\ref{fig:resd_conv} illustrates the convergence of the corresponding greedy algorithm which stops within 10 iterations. The effectivity in Fig.~\ref{fig:resd_eff} is of order $\mathcal{O}(10)$. The algorithm results in a ROM of size 20.\\
\subsubsection{Test 3. Randomized error estimator}
To compute the randomized estimator $\tilde \Delta(\mu)$ from~\cite{morSemZP18}, $K=20$ random vectors are generated from a normal distribution by the MATLAB\textsuperscript \textregistered command \texttt{mvnrnd} with seed $0$, using the \texttt{MersenneTwister} random number generator. The projection matrix $V_{rd}$ is constructed for the $K$ random dual systems in~(\ref{eq:dual_rand}), using Algorithm 3.1 in~\cite{morSemZP18} with $\texttt{tol}_{rd} = 0.5$. It is seen from Fig.~\ref{fig:rand_conv}  that the greedy algorithm takes $8$ iterations to converge leading to a ROM of order 16.  In Fig.~\ref{fig:rand_eff}, the effectivity of $\tilde \Delta(\mu)$ is plotted. Although it is close to one, there are obvious underestimations at some iterations. Furthermore, generating $V_{rd}$ takes considerable time as shown in Table \ref{tab:offtime}, even for relatively large error tolerances ($\thicksim\mathcal{O}(10^{-1})$).\\
\subsubsection{Test 4. Proposed error estimator}
\begin{figure}[tbp]
	\centering
	\input{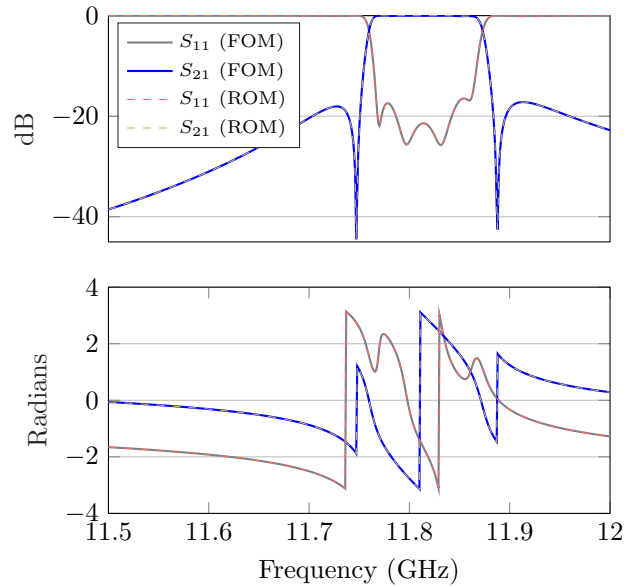}
	\caption{Dual-mode circular waveguide filter: scattering parameter responses of FOM and ROM computed from Test 4. Top: The magnitudes of the scattering parameter responses $|S_{11}|$, $|S_{21}|$; Bottom: The phases $\angle(S_{11}), \angle(S_{21})$.}
	\label{fig:DualModeFilter_Sparameters}
\end{figure}

The performance of $\|\tilde{e}(s)\|_2$ is tested by using Algorithm 1. As illustrated in Fig.~\ref{fig:DualModeFilter_Test4}, Algorithm 1 needs $8$ iterations to converge with effectivity very close to $1$. Compared with the standard error estimator $\delta(\mu)$ and the residual norm error estimator $\|r(\mu)\|_2$, $\|\tilde{e}(s)\|_2$ is much tighter, and therefore converges faster leading to a ROM of smaller order $r=16$. It is also more reliable than the randomized estimator $\tilde \Delta(\mu)$ with almost no underestimation. Fig.~\ref{fig:DualModeFilter_Sparameters} provides a comparison of the scattering parameter response for this filter resulting from the FOM and the ROM. As can be seen, the ROM obtained using Test 4 produces a very good match.

\begin{figure}[tbp]
	\centering
	\subfloat[]{\label{fig:std_conv}
%
%
\definecolor{mycolor1}{rgb}{0.00000,0.44700,0.74100}%
\begin{tikzpicture}

\begin{axis}[%
width=\fwidth,
height=\fheight,
at={(0\fwidth,0\fheight)},
scale only axis,
xmin=1,
xmax=10,
xlabel style={font=\color{white!15!black}},
xlabel={Iterations},
ymode=log,
ymin=1e-10,
ymax=10000,
yminorticks=true,
ytick = {10000, 0.1, 0.000001},
ylabel style={font=\color{white!15!black}},
ylabel={Maximum error},
axis background/.style={fill=white},
xmajorgrids,
ymajorgrids,
yminorgrids,
legend style={legend cell align=left, align=left, draw=white!15!black}
]
\addplot [color=gray, thick]
  table[row sep=crcr]{%
1	3139.0785178615\\
2	4346.95438448382\\
3	1576.20218781667\\
4	5.08871394920877\\
5	0.0255573253112913\\
6	0.748981714612039\\
7	0.000649701317834131\\
8	3.36436383914769e-05\\
9	2.17306559490696e-05\\
10	1.69369050053765e-07\\
};
\addlegendentry{$\epsilon_{\text{est}}$}

\addplot [color=black, dashed, mark=diamond, mark options={solid, black}, thick]
  table[row sep=crcr]{%
1	185.786607643769\\
2	75.2527407224027\\
3	3.43334353085329\\
4	0.0129105637589105\\
5	0.00166568448421715\\
6	0.0301848153302552\\
7	2.56764575745013e-05\\
8	2.76568306206683e-07\\
9	8.52177584450135e-07\\
10	2.72492787490513e-09\\
};
\addlegendentry{$\epsilon_{\text{true}}$}

\end{axis}
\end{tikzpicture}%
} 
	\subfloat[]{\label{fig:std_eff}
%
%
\definecolor{mycolor1}{rgb}{0.00000,0.44700,0.74100}%
\begin{tikzpicture}

\begin{axis}[%
width=\fwidth,
height=\fheight,
at={(0\fwidth,0\fheight)},
scale only axis,
xmin=1,
xmax=10,
xlabel style={font=\color{white!15!black}},
xlabel={Iterations},
ymin=0,
ymax=500,
ylabel style={font=\color{white!15!black}},
ylabel={\texttt{eff}},
axis background/.style={fill=white},
xmajorgrids,
ymajorgrids
]
\addplot [color=white!20!black, thick]
  table[row sep=crcr]{%
1	16.8961506842325\\
2	57.764731792549\\
3	459.086652312048\\
4	394.15118071019\\
5	15.3434372196262\\
6	24.8131951916005\\
7	25.3033860278037\\
8	121.646760082244\\
9	25.5001496702019\\
10	62.1554249613529\\
};
\end{axis}
\end{tikzpicture}%
}
	\caption{Dual-mode circular waveguide filter: results for Test~1. (a) Convergence of the greedy algorithm. (b) Effectivity (\texttt{eff}).}
	\label{fig:DualModeFilter_Test1}
\end{figure}
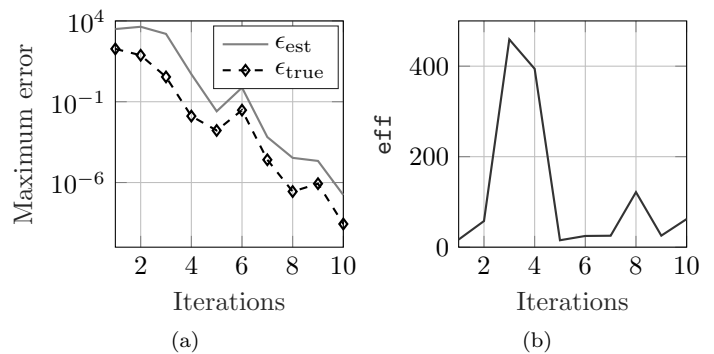	
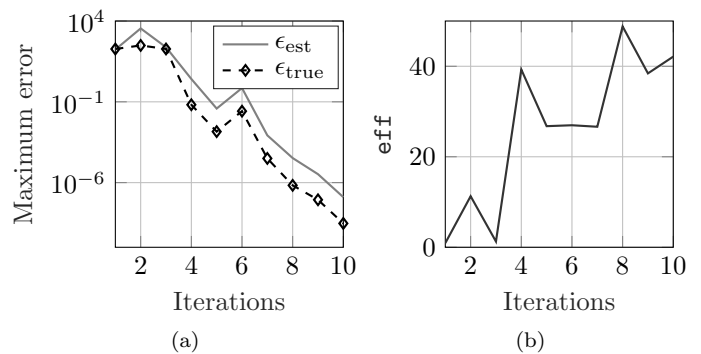
\begin{figure}[tbp]
	\centering
	\subfloat[]{\label{fig:resd_conv}
%
%
\definecolor{mycolor1}{rgb}{0.00000,0.44700,0.74100}%
\begin{tikzpicture}

\begin{axis}[%
width=\fwidth,
height=\fheight,
at={(0\fwidth,0\fheight)},
scale only axis,
xmin=1,
xmax=10,
xlabel style={font=\color{white!15!black}},
xlabel={Iterations},
ymode=log,
ymin=1e-10,
ymax=10000,
yminorticks=true,
ytick = {10000, 0.1, 0.000001},
ylabel style={font=\color{white!15!black}},
ylabel={Maximum error},
axis background/.style={fill=white},
xmajorgrids,
ymajorgrids,
yminorgrids,
legend style={legend cell align=left, align=left, draw=white!15!black}
]
\addplot [color=gray, thick]
  table[row sep=crcr]{%
1	171.572242499304\\
2	3482.00244618299\\
3	238.318837370104\\
4	2.59670641211638\\
5	0.0384244181842813\\
6	0.71166784919471\\
7	0.000839821625520888\\
8	3.32321537222735e-05\\
9	3.33633459643698e-06\\
10	1.26368897099303e-07\\
};
\addlegendentry{$\epsilon_{\text{est}}$}

\addplot [color=black, dashed, mark=diamond, mark options={solid, black},thick]
  table[row sep=crcr]{%
1	185.786607643769\\
2	308.592848067113\\
3	186.381087450447\\
4	0.0660998098874641\\
5	0.00143571422223379\\
6	0.026390516085957\\
7	3.15320989259145e-05\\
8	6.80973799542636e-07\\
9	8.68222090877622e-08\\
10	2.99887664671263e-09\\
};
\addlegendentry{$\epsilon_{\text{true}}$}

\end{axis}
\end{tikzpicture}%
} 
	\subfloat[]{\label{fig:resd_eff}
%
%
\definecolor{mycolor1}{rgb}{0.00000,0.44700,0.74100}%
\begin{tikzpicture}

\begin{axis}[%
width=\fwidth,
height=\fheight,
at={(0\fwidth,0\fheight)},
scale only axis,
xmin=1,
xmax=10,
xlabel style={font=\color{white!15!black}},
xlabel={Iterations},
ymin=0,
ymax=50,
ylabel style={font=\color{white!15!black}},
ylabel={\texttt{eff}},
axis background/.style={fill=white},
xmajorgrids,
ymajorgrids
]
\addplot [color=white!20!black, thick]
  table[row sep=crcr]{%
1	0.923490905373979\\
2	11.2834839433016\\
3	1.27866427130631\\
4	39.2846275433668\\
5	26.7632775306062\\
6	26.9668030316923\\
7	26.6338637175429\\
8	48.8009285299864\\
9	38.4272023424851\\
10	42.1387445988578\\
};
\end{axis}
\end{tikzpicture}%
}
	\caption{Dual-mode circular waveguide filter: results for Test~2. (a) Convergence of the greedy algorithm. (b) Effectivity (\texttt{eff}).}
	\label{fig:DualModeFilter_Test2}
\end{figure}
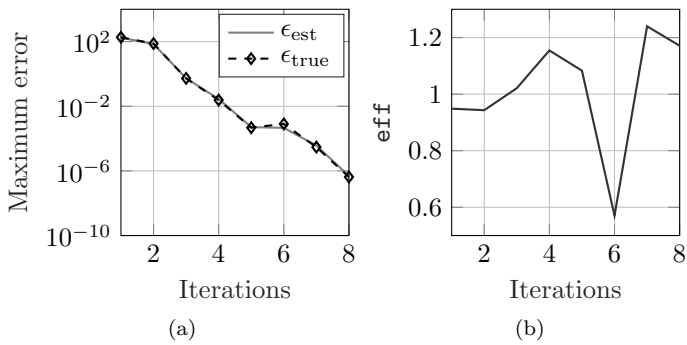
\begin{figure}[tbp]
	\centering
	\subfloat[]{\label{fig:rand_conv}
%
%
\definecolor{mycolor1}{rgb}{0.00000,0.44700,0.74100}%
\begin{tikzpicture}

\begin{axis}[%
width=\fwidth,
height=\fheight,
at={(0\fwidth,0\fheight)},
scale only axis,
xmin=1,
xmax=8,
xlabel style={font=\color{white!15!black}},
xlabel={Iterations},
ymode=log,
ymin=1e-10,
ymax=10000,
ytick = {0.0000000001, 0.000001, 0.01, 100},
yminorticks=true,
ylabel style={font=\color{white!15!black}},
ylabel={Maximum error},
axis background/.style={fill=white},
xmajorgrids,
ymajorgrids,
yminorgrids,
legend style={legend cell align=left, align=left, draw=white!15!black}
]
\addplot [color=gray,thick]
  table[row sep=crcr]{%
1	176.243554676942\\
2	70.9810089696426\\
3	0.540615427193463\\
4	0.0278822395747491\\
5	0.000510715714635458\\
6	0.000460375867747043\\
7	3.68617137637675e-05\\
8	4.9114233807274e-07\\
};
\addlegendentry{$\epsilon_{\text{est}}$}

\addplot [color=black, dashed, mark=diamond, mark options={solid, black},thick]
  table[row sep=crcr]{%
1	185.786607643769\\
2	75.2527407224027\\
3	0.529485895700952\\
4	0.0241512671774332\\
5	0.000471581640574201\\
6	0.000806816700981526\\
7	2.97323381223845e-05\\
8	4.19370735203354e-07\\
};
\addlegendentry{$\epsilon_{\text{true}}$}

\end{axis}
\end{tikzpicture}%
} 
	\subfloat[]{\label{fig:rand_eff}
%
%
\definecolor{mycolor1}{rgb}{0.00000,0.44700,0.74100}%
\begin{tikzpicture}

\begin{axis}[%
width=\fwidth,
height=\fheight,
at={(0\fwidth,0\fheight)},
scale only axis,
xmin=1,
xmax=8,
xlabel style={font=\color{white!15!black}},
xlabel={Iterations},
ymin=0.5,
ymax=1.3,
ylabel style={font=\color{white!15!black}},
ylabel={\texttt{eff}},
axis background/.style={fill=white},
xmajorgrids,
ymajorgrids
]
\addplot [color=white!20!black, thick]
  table[row sep=crcr]{%
1	0.948634333293146\\
2	0.94323486810244\\
3	1.02101950511406\\
4	1.15448350473312\\
5	1.08298472776338\\
6	0.57060775661557\\
7	1.2397852335742\\
8	1.1711411809281\\
};
\end{axis}
\end{tikzpicture}%
}
	\caption{Dual-mode circular waveguide filter: results for Test~3. (a) Convergence of the greedy algorithm. (b) Effectivity (\texttt{eff}).}
	\label{fig:DualModeFilter_Test3}
\end{figure} 
\begin{figure}[t!]
	\centering
	\subfloat[]{\label{fig:est_conv}
%
%
\definecolor{mycolor1}{rgb}{0.00000,0.44700,0.74100}%
\begin{tikzpicture}

\begin{axis}[%
width=\fwidth,
height=\fheight,
at={(0\fwidth,0\fheight)},
scale only axis,
xmin=1,
xmax=8,
xlabel style={font=\color{white!15!black}},
xlabel={Iterations},
ymode=log,
ymin=1e-10,
ymax=10000,
ytick = {0.0000000001, 0.000001, 0.01, 100},
yminorticks=true,
ylabel style={font=\color{white!15!black}},
ylabel={Maximum error},
axis background/.style={fill=white},
xmajorgrids,
ymajorgrids,
yminorgrids,
legend style={legend cell align=left, align=left, draw=white!15!black}
]
\addplot [color=gray,thick]
  table[row sep=crcr]{%
1	957.331377929503\\
2	214.036926537715\\
3	987.354649417546\\
4	0.042762443444301\\
5	0.000642802722605922\\
6	0.000357854041576229\\
7	3.17866489071252e-06\\
8	3.80015613436873e-07\\
};
\addlegendentry{$\epsilon_{\text{est}}$}

\addplot [color=black, dashed, mark=diamond, mark options={solid, black},thick]
  table[row sep=crcr]{%
1	185.786607643769\\
2	214.036885408335\\
3	987.354649418769\\
4	0.0427624435057797\\
5	0.000642802722607191\\
6	0.000362367914551706\\
7	3.17866483739295e-06\\
8	3.80015672049463e-07\\
};
\addlegendentry{$\epsilon_{\text{true}}$}

\end{axis}
\end{tikzpicture}%
} 
	\subfloat[]{\label{fig:est_eff}
%
%
\definecolor{mycolor1}{rgb}{0.00000,0.44700,0.74100}%
\begin{tikzpicture}

\begin{axis}[%
width=\fwidth,
height=\fheight,
at={(0\fwidth,0\fheight)},
scale only axis,
xmin=1,
xmax=8,
xlabel style={font=\color{white!15!black}},
xlabel={Iterations},
ymin=0.8,
ymax=5,
ylabel style={font=\color{white!15!black}},
ylabel={\texttt{eff}},
axis background/.style={fill=white},
xmajorgrids,
ymajorgrids
]
\addplot [color=white!20!black, thick]
  table[row sep=crcr]{%
1	5.15285461137819\\
2	1.00000019216025\\
3	0.999999999998761\\
4	0.999999998562321\\
5	0.999999999998025\\
6	0.987543397761742\\
7	1.0000000167742\\
8	0.999999845762705\\
};
\end{axis}
\end{tikzpicture}%
}
	\caption{Dual-mode circular waveguide filter: results for Test~4. (a) Convergence of Algorithm~\ref{alg:rom_errest}. (b) Effectivity (\texttt{eff}).}
	\label{fig:DualModeFilter_Test4}
\end{figure} 
%
%
For the next set of tests, we consider three different antenna models: 
\begin{enumerate}[(i)]
	\item Substrate Integrated Waveguide (SIW) antenna.
	\item Antipodal Vivaldi (AV) antenna.
	\item Dielectric Resonator (DR) antenna.
\end{enumerate}
For each of the models, we perform Tests 2, 3 and 4, as for Example I. For these examples, performing Test 1 is computationally expensive; moreover, the resulting estimation of the error is not sharp, just like in Example I. Therefore, for the remaining examples, we do not show the additional results involving the standard error estimator.

\subsection{Example II : Substrate Integrated Waveguide Antenna}
\begin{figure}[tbp]
	\centering
	\includegraphics[width=\linewidth]{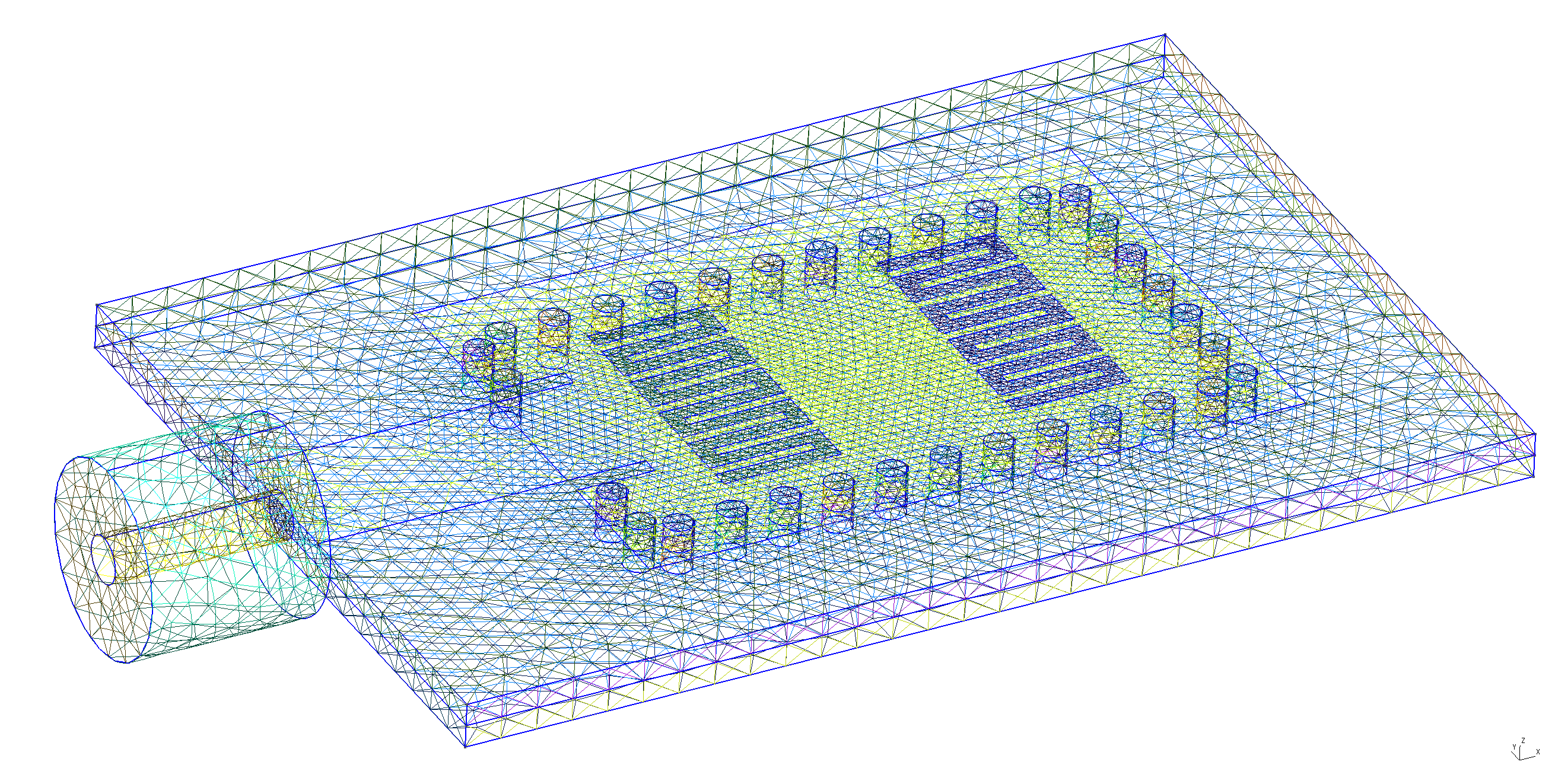}
	\caption{Geometry of the SIW antenna proposed in \cite{DonI10}.}
	\label{fig:SIW_Geometry}
\end{figure}

Such antennas have gained popularity recently owing to their low-cost and efficiency \cite{DonI10}. Fig.~\ref{fig:SIW_Geometry} shows the model of the antenna along with the mesh used for its discretization. 
The system is of order $n = 390,302$, with just one input. The frequency band of interest is $[6, \, 9 ]$ GHz. The training set $\Xi$ consists of $61$ uniformly sampled parameters from this band. The tolerance ($\texttt{tol}$) for the greedy algorithm is set as $10^{-4}$.\\
\subsubsection{Test 2. Residual norm as a heuristic error estimator}
We illustrate the convergence of the greedy algorithm using the residual estimator in Fig.~\ref{fig:SIW_resd_conv}. The algorithm converges in 8 iterations to a ROM of size $r = 16$ and takes around $12$ minutes. The effectivity shown in Fig.~\ref{fig:SIW_resd_eff} is less than one for the first four iterations, meaning the estimator underestimates the true error. For the next four iterations, it is of order $\mathcal{O}(10)$. Using the residual norm as an error estimator enjoys the advantage that it is very easy to implement, with only marginal overhead costs in terms of computation. However, in general, it is a crude approximation to the actual error.\\
\subsubsection{Test 3. Randomized error estimator}
Next, we use the error estimator from \cite{morSemZP18}. We use $K = 5$ randomly distributed vectors, each of length $n$ to construct $V_{rd}$. The random number generator \texttt{MersenneTwister} was used with the seed set to $1$. The tolerance used to obtain $V_{rd}$ is $\texttt{tol}_{rd} = 1$. Even for such a crude tolerance, the algorithm requires around 28 minutes to converge. The long time to converge and the lack of any \textit{a priori} knowledge to choose the number of random vectors $K$ and the tolerance are the main disadvantages of using this algorithm.  Fig.~\ref{fig:SIW_Test3} shows the convergence of the greedy algorithm and the effectivity. The overall procedure requires $7$ iterations to converge to a ROM of dimension $r = 14$.\\
\subsubsection{Test 4. Proposed error estimator}
For the proposed error estimator, we show the results of Algorithm 1 in Fig.~\ref{fig:SIW_Test4}. It requires $7$ iterations to converge to a ROM with $r = 14$. The overall time to converge was around 21 minutes, much less than the time required to precompute $V_{rd}$ in Test 3. Further, compared to the other methods, the proposed method achieves a better effectivity, where only at the first two iterations, the error estimator underestimates the true error. This is quite reasonable, since at those stages the ROM of the residual system (\ref{eq:sys_resi}) is not yet accurate enough, leading to large $\gamma(\mu)$ in (\ref{eq:err_bound}). In Fig.~\ref{fig:SIW_Sparameters} we show the scattering parameter response at the coaxial port obtained from the FOM and the ROM using Test 4. The results are nearly identical, thus showing the good approximation capabilities of the proposed approach.
\begin{figure}[tbp]
	\centering
	\subfloat[]{\label{fig:SIW_resd_conv}
%
%
\definecolor{mycolor1}{rgb}{0.00000,0.44700,0.74100}%
\begin{tikzpicture}

\begin{axis}[%
width=\fwidth,
height=\fheight,
at={(0\fwidth,0\fheight)},
scale only axis,
xmin=1,
xmax=8,
xlabel style={font=\color{white!15!black}},
xlabel={Iterations},
ymode=log,
ymin=1e-06,
ymax=10000,
yminorticks=true,
ytick = {100, 0.1, 0.0001},
ylabel style={font=\color{white!15!black}},
ylabel={Maximum error},
axis background/.style={fill=white},
xmajorgrids,
ymajorgrids,
yminorgrids,
legend style={legend cell align=left, align=left, draw=white!15!black}
]
\addplot [color=gray,thick]
  table[row sep=crcr]{%
1	31.5295418394526\\
2	11.9813863649667\\
3	30.4990479246111\\
4	18.1909752866828\\
5	1.79153522076753\\
6	0.037168636175321\\
7	0.00102507376144201\\
8	7.07300989274003e-05\\
};
\addlegendentry{$\epsilon_{\text{est}}$}

\addplot [color=black, dashed, mark=diamond, mark options={solid, black},thick]
  table[row sep=crcr]{%
1	47.5414818924125\\
2	47.4656301476196\\
3	47.4777229902058\\
4	44.4987227215323\\
5	0.0765206719591924\\
6	0.00116256526450628\\
7	3.40695381519223e-05\\
8	2.25673586318337e-06\\
};
\addlegendentry{$\epsilon_{\text{true}}$}
\end{axis}
\end{tikzpicture}%
} 
	\subfloat[]{\label{fig:SIW_resd_eff}
%
%
\definecolor{mycolor1}{rgb}{0.00000,0.44700,0.74100}%
\begin{tikzpicture}

\begin{axis}[%
width=\fwidth,
height=\fheight,
at={(0\fwidth,0\fheight)},
scale only axis,
xmin=1,
xmax=8,
xlabel style={font=\color{white!15!black}},
xlabel={Iterations},
ymin=0,
ymax=35,
ylabel style={font=\color{white!15!black}},
ylabel={\texttt{eff}},
axis background/.style={fill=white},
xmajorgrids,
ymajorgrids
]
\addplot [color=white!20!black, thick]
  table[row sep=crcr]{%
1	0.663200653080287\\
2	0.252422359667494\\
3	0.642386492100784\\
4	0.408797695172506\\
5	23.4124345081933\\
6	31.9712254529691\\
7	30.0876917342119\\
8	31.3417711311716\\
};
\end{axis}
\end{tikzpicture}%
}
	\caption{SIW antenna: results for Test 2. (a) Convergence of the greedy algorithm. (b) Effectivity (\texttt{eff}).}
	\label{fig:SIW_Test2}
\end{figure}
\begin{figure}[tbp]
	\centering
	\subfloat[]{\label{fig:SIW_rand_conv}
%
%
\definecolor{mycolor1}{rgb}{0.00000,0.44700,0.74100}%
\begin{tikzpicture}

\begin{axis}[%
width=\fwidth,
height=\fheight,
at={(0\fwidth,0\fheight)},
scale only axis,
xmin=1,
xmax=7,
xtick = {1,3,5,7},
xlabel style={font=\color{white!15!black}},
xlabel={Iterations},
ymode=log,
ymin=1e-06,
ymax=10000,
yminorticks=true,
ytick = {100, 0.1, 0.0001},
ylabel style={font=\color{white!15!black}},
ylabel={Maximum error},
axis background/.style={fill=white},
xmajorgrids,
ymajorgrids,
yminorgrids,
legend style={legend cell align=left, align=left, draw=white!15!black}
]
\addplot [color=gray,thick]
  table[row sep=crcr]{%
1	16.555004579316\\
2	13.3870818749191\\
3	8.51364011435234\\
4	0.376698204907059\\
5	0.0173178571089941\\
6	0.000523461802802976\\
7	2.55583337996711e-05\\
};
\addlegendentry{$\epsilon_{\text{est}}$}

\addplot [color=black, dashed, mark=diamond, mark options={solid, black},thick]
  table[row sep=crcr]{%
1	47.5414818924125\\
2	21.7550229210538\\
3	7.58187623714402\\
4	0.611723156174596\\
5	0.0195533604448549\\
6	0.000569414255062329\\
7	3.32960409715024e-05\\
};
\addlegendentry{$\epsilon_{\text{true}}$}

\end{axis}
\end{tikzpicture}%
} 
	\subfloat[]{\label{fig:SIW_rand_eff}
%
%
\definecolor{mycolor1}{rgb}{0.00000,0.44700,0.74100}%
\begin{tikzpicture}

\begin{axis}[%
width=\fwidth,
height=\fheight,
at={(0\fwidth,0\fheight)},
scale only axis,
xmin=1,
xmax=7,
xtick = {1,3,5,7},
xlabel style={font=\color{white!15!black}},
xlabel={Iterations},
ymin=0.3,
ymax=1.2,
ylabel style={font=\color{white!15!black}},
ylabel={\texttt{eff}},
axis background/.style={fill=white},
xmajorgrids,
ymajorgrids
]
\addplot [color=white!20!black, thick]
  table[row sep=crcr]{%
1	0.348222308609992\\
2	0.615355907621846\\
3	1.12289357516066\\
4	0.615798504772546\\
5	0.885671655152809\\
6	0.919298732248418\\
7	0.767608792334983\\
};
\end{axis}
\end{tikzpicture}%
}
	\caption{SIW antenna: results for Test 3. (a) Convergence of the greedy algorithm. (b) Effectivity (\texttt{eff}).}
	\label{fig:SIW_Test3}
\end{figure}
\begin{figure}[tbp]
	\centering
	\subfloat[]{\label{fig:SIW_prop_conv}
%
%
\definecolor{mycolor1}{rgb}{0.00000,0.44700,0.74100}%
\begin{tikzpicture}

\begin{axis}[%
width=\fwidth,
height=\fheight,
at={(0\fwidth,0\fheight)},
scale only axis,
xmin=1,
xmax=7,
xtick = {1,3,5,7},
xlabel style={font=\color{white!15!black}},
xlabel={Iterations},
ymode=log,
ymin=1e-06,
ymax=10000,
yminorticks=true,
ytick = {100, 0.1, 0.0001},
ylabel style={font=\color{white!15!black}},
ylabel={Maximum error},
axis background/.style={fill=white},
xmajorgrids,
ymajorgrids,
yminorgrids,
legend style={legend cell align=left, align=left, draw=white!15!black}
]
\addplot [color=gray, thick]
  table[row sep=crcr]{%
1	7.56099080344208\\
2	25.7391503513823\\
3	47.8509976549763\\
4	0.611723229018772\\
5	0.0155526159138844\\
6	0.00078205796941349\\
7	2.22839034977085e-05\\
};
\addlegendentry{$\epsilon_{\text{est}}$}
\addplot [color=black, dashed, mark=diamond, mark options={solid, black},thick]
  table[row sep=crcr]{%
1	47.5414818924125\\
2	47.4656301476196\\
3	47.2069111086789\\
4	0.61172315617875\\
5	0.0155526160424621\\
6	0.000782057974877403\\
7	2.2283903857221e-05\\
};
\addlegendentry{$\epsilon_{\text{true}}$}
\end{axis}
\end{tikzpicture}%
} 
	\subfloat[]{\label{fig:SIW_prop_eff}
%
%
\definecolor{mycolor1}{rgb}{0.00000,0.44700,0.74100}%
\begin{tikzpicture}

\begin{axis}[%
width=0.951\fwidth,
height=\fheight,
at={(0\fwidth,0\fheight)},
scale only axis,
xmin=1,
xmax=7,
xtick = {1,3,5,7},
xlabel style={font=\color{white!15!black},font=\large},
xlabel={Iterations},
ymin=0.05,
ymax=1.02,
ylabel style={font=\color{white!15!black},font=\large},
ylabel={\texttt{eff}},
axis background/.style={fill=white},
xmajorgrids,
ymajorgrids
]
\addplot [color=white!20!black, thick]
  table[row sep=crcr]{%
1	0.159039863766821\\
2	0.542269222410674\\
3	1.01364390363976\\
4	1.00000011907351\\
5	0.999999991732729\\
6	0.999999993013416\\
7	0.999999983866721\\
};
\end{axis}
\end{tikzpicture}%
}
	\caption{SIW antenna: results for Test 4. (a) Convergence of Algorithm~\ref{alg:rom_errest}. (b) Effectivity (\texttt{eff}).}
	\label{fig:SIW_Test4}
\end{figure}
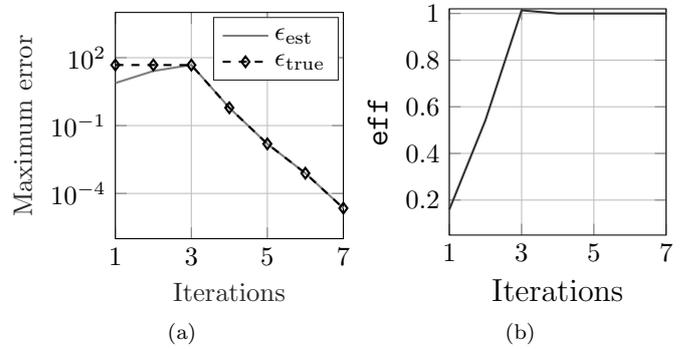
\begin{figure}[tbp]
	\centering
	\input{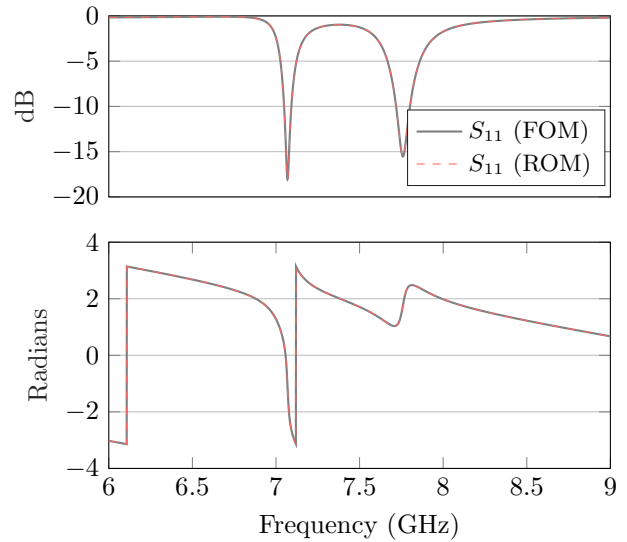}
	\caption{SIW antenna: scattering parameter responses of FOM and ROM computed from Test 4. Top: The magnitude of the scattering parameter response $|S_{11}|$; Bottom: The phase $\angle(S_{11})$.}
	\label{fig:SIW_Sparameters}
\end{figure}
\subsection{Example III : Antipodal Vivaldi Antenna}
For the next example, we employ an Antipodal Vivaldi antenna. It is known for good wide-band impedance performance \cite{LouJ05} and contains a large number of in-band resonances. The discretized model, shown in Fig.~\ref{fig:AVA_Geometry} is of dimension $n = 283,846$ with $f \in [1, \, 6]$ GHz being the parameter space of interest. The tolerance for the ROM is set as $10^{-3}$. The training set is made up of $51$ uniform samples from the parameter range of interest.\\
\subsubsection{Test 2. Residual norm as a heuristic error estimator}
In Fig.~\ref{fig:AVA_resd_conv} the convergence of the greedy algorithm is shown. The algorithm requires up to $20$ iterations to achieve the set tolerance, taking roughly $30$ minutes. For this model, the norm of the residual overestimates the true error by roughly one order of magnitude, as seen in Fig.~\ref{fig:AVA_resd_eff}, where the effectivity is illustrated. The resulting ROM has dimension $r = 40$.\\
\subsubsection{Test 3. Randomized error estimator}
To construct $V_{rd}$ for this example, we draw $K = 6$ random vectors using the \texttt{mvnrnd} command and set $\texttt{tol}_{rd}$ to be $1$. The random number generator \texttt{MersenneTwister} was used with the seed set to $1$. The procedure based on Algorithm 3.1 from \cite{morSemZP18} to generate $V_{rd}$ takes nearly $89$ minutes.
The results of the greedy algorithm are displayed in Fig.~\ref{fig:AVA_Test3}. Compared to Test 2 for this example, Test 3 requires only $17$ iterations to converge. It results in a ROM of dimension $r = 34$.\\
\subsubsection{Test 4. Proposed error estimator}
Using Algorithm 1 with \texttt{tol} = $10^{-3}$ for the Antipodal Vivaldi antenna results in a ROM of dimension $r = 36$, with the convergence achieved in $18$ iterations. The convergence of the maximum error and the corresponding effectivity are shown in Fig.~\ref{fig:AVA_Test4}. Although in comparison to the randomized error estimator in Test 3, the proposed approach takes one extra iteration to converge, the overall time for Test 4 is only $56$ minutes. Moreover, the effectivity is nearly $1$ for most of the iterations. The randomized estimator tending to underestimate the true error also explains why it takes fewer iterations to converge: the algorithm stops according to the already small error estimate, even before the true error is below the tolerance over the whole frequency domain.

While we have set a coarser tolerance of $10^{-3}$ for this example, this is not a limitation of Algorithm 1. Depending on the application or problem at hand, there is usually a trade-off between accuracy and computational cost. The tolerance can be used to define this trade-off. In some applications, a greater accuracy for the ROM may be desired in order to perform highly accurate design choices by simulating the ROM. In such scenarios, the tolerance can be set low to achieve a high-quality ROM. However, there also arise situations where a highly accurate ROM is not needed or the offline computational cost that can be afforded is less. In such cases, a higher tolerance can be set to achieve the desired trade-off. 

For this example, we run Test 4 by setting a lower tolerance of $10^{-4}$. As can be seen from Fig.~\ref{fig:AVA_prop_conv_smalltol}, the greedy algorithm converges to this smaller tolerance in $20$ iterations, taking $62$ minutes. We also note that the effectivity shown in Fig.~\ref{fig:AVA_prop_eff_smalltol} is also nearly identical and close to $1$ for most of the iterations. Finally, we make a comparison of the scattering parameter response of the Vivaldi antenna obtained from the FOM and the ROM. The results in Fig.~\ref{fig:AVA_Sparameters} display a very close match between the two.
\begin{figure}[tbp]
	\centering
	\includegraphics[width=\linewidth]{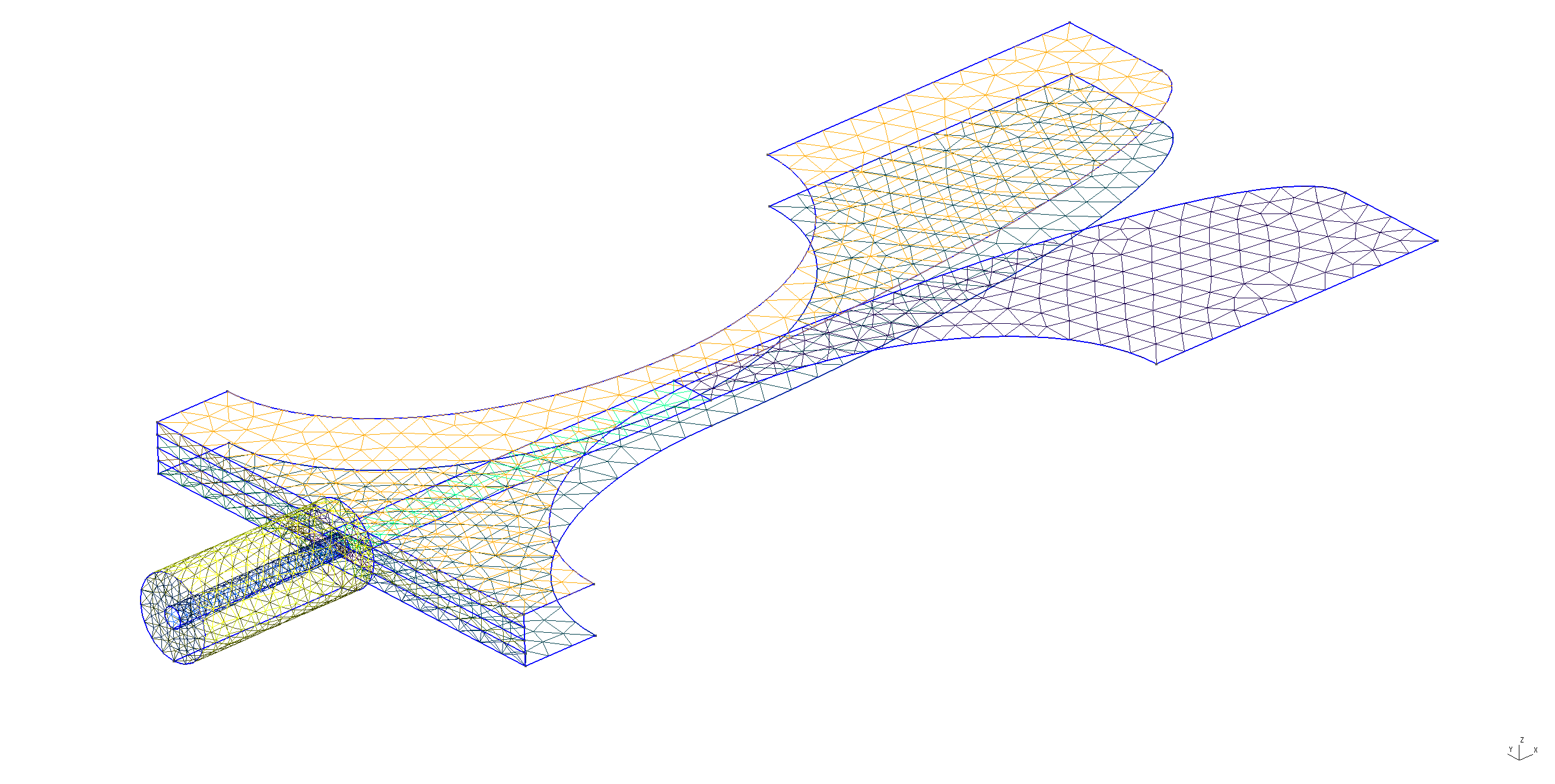} 
	\caption{Antipodal Vivaldi antenna detailed in \cite{LouJ05}.}
	\label{fig:AVA_Geometry}
\end{figure}
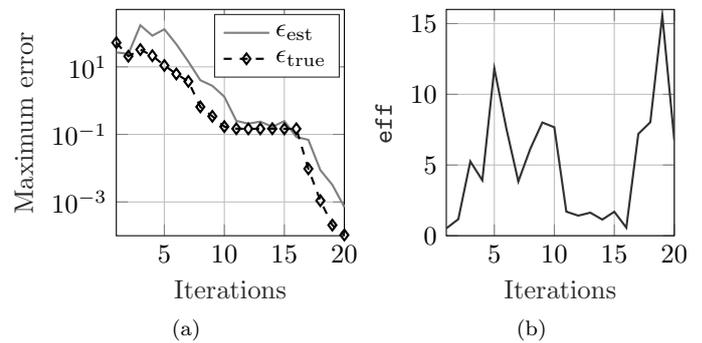
\begin{figure}[tbp]
	\centering
	\subfloat[]{\label{fig:AVA_resd_conv}
%
%
\definecolor{mycolor1}{rgb}{0.00000,0.44700,0.74100}%
\begin{tikzpicture}

\begin{axis}[%
width=\fwidth,
height=\fheight,
at={(0\fwidth,0\fheight)},
scale only axis,
xmin=1,
xmax=20,
xlabel style={font=\color{white!15!black}},
xlabel={Iterations},
ymode=log,
ymin=0.0001,
ymax=500,
yminorticks=true,
ytick = {10,0.1,0.001},
ylabel style={font=\color{white!15!black}},
ylabel={Maximum error},
axis background/.style={fill=white},
xmajorgrids,
ymajorgrids,
yminorgrids,
legend style={legend cell align=left, align=left, draw=white!15!black}
]
\addplot [color=gray,thick]
  table[row sep=crcr]{%
1	26.8469799767679\\
2	24.0904784784924\\
3	170.057009648333\\
4	83.6960869767928\\
5	129.753136272033\\
6	46.7229749412029\\
7	14.2920343396976\\
8	4.01377477929043\\
9	2.74583930615883\\
10	1.30271180043834\\
11	0.251234900394967\\
12	0.20747576103024\\
13	0.23900828630929\\
14	0.168081783939538\\
15	0.247542415331774\\
16	0.0838908292265715\\
17	0.0690322883150573\\
18	0.00876068919633795\\
19	0.00320974897671287\\
20	0.000708780045384921\\
};
\addlegendentry{$\epsilon_{\text{est}}$}

\addplot [color=black, dashed, mark=diamond, mark options={solid,black},thick]
  table[row sep=crcr]{%
1	51.8466903471054\\
2	20.5987191782988\\
3	32.3685625796358\\
4	21.3217423867864\\
5	10.9731861653704\\
6	6.15220212294054\\
7	3.71865141646668\\
8	0.652498469925707\\
9	0.342947278342973\\
10	0.169855286378299\\
11	0.147652450261356\\
12	0.146203063027784\\
13	0.146522771816691\\
14	0.147339304806147\\
15	0.145870278670295\\
16	0.146573606033582\\
17	0.00959088230402838\\
18	0.00109166405896595\\
19	0.000206240303919245\\
20	0.000105137980370237\\
};
\addlegendentry{$\epsilon_{\text{true}}$}

\end{axis}
\end{tikzpicture}%
} 
	\subfloat[]{\label{fig:AVA_resd_eff}
%
%
\definecolor{mycolor1}{rgb}{0.00000,0.44700,0.74100}%
\begin{tikzpicture}

\begin{axis}[%
width=\fwidth,
height=\fheight,
at={(0\fwidth,0\fheight)},
scale only axis,
xmin=1,
xmax=20,
xlabel style={font=\color{white!15!black}},
xlabel={Iterations},
ymin=0,
ymax=16,
ylabel style={font=\color{white!15!black}},
ylabel={\texttt{eff}},
axis background/.style={fill=white},
xmajorgrids,
ymajorgrids
]
\addplot [color=white!20!black, thick]
  table[row sep=crcr]{%
1	0.517814730256292\\
2	1.16951341828439\\
3	5.25377082253637\\
4	3.92538684027349\\
5	11.824563469224\\
6	7.59451234005799\\
7	3.84333801130448\\
8	6.15139339674994\\
9	8.00659308167124\\
10	7.66953933677967\\
11	1.70152882630977\\
12	1.41909312112573\\
13	1.6312023267503\\
14	1.14078035158834\\
15	1.69700378712023\\
16	0.572346082604742\\
17	7.19769945316316\\
18	8.02507797557817\\
19	15.563150924999\\
20	6.74142724531129\\
};
\end{axis}
\end{tikzpicture}%
}
	\caption{Antipodal Vivaldi antenna: results for Test 2. (a) Convergence of the greedy algorithm. (b) Effectivity (\texttt{eff}).}
	\label{fig:AVA_Test2}
\end{figure}
\begin{figure}[tbp]
	\centering
	\subfloat[]{\label{fig:AVA_rand_conv}
%
%
\definecolor{mycolor1}{rgb}{0.00000,0.44700,0.74100}%
\begin{tikzpicture}

\begin{axis}[%
width=\fwidth,
height=\fheight,
at={(0\fwidth,0\fheight)},
scale only axis,
xmin=1,
xmax=17,
xlabel style={font=\color{white!15!black}},
xlabel={Iterations},
ymode=log,
ymin=0.0001,
ymax=100,
yminorticks=true,
ytick = {10,0.1,0.001},
ylabel style={font=\color{white!15!black}},
ylabel={Maximum error},
axis background/.style={fill=white},
xmajorgrids,
ymajorgrids,
yminorgrids,
legend style={legend cell align=left, align=left, draw=white!15!black}
]
\addplot [color=gray,thick]
  table[row sep=crcr]{%
1	36.7928225705684\\
2	21.7663572106117\\
3	10.6541260762615\\
4	9.50802601338422\\
5	8.31988506332221\\
6	1.43802223851536\\
7	1.15572009951063\\
8	0.341552645668893\\
9	0.172461915157393\\
10	0.4407623185243\\
11	0.291790675809762\\
12	0.121561805430603\\
13	0.0546092203324502\\
14	0.0567314754954377\\
15	0.0162739415339935\\
16	0.00610351965648817\\
17	0.000993884219526163\\
};
\addlegendentry{$\epsilon_{\text{est}}$}

\addplot [color=black, dashed, mark=diamond, mark options={solid, black},thick]
  table[row sep=crcr]{%
1	51.8466903471054\\
2	20.5832866911299\\
3	10.6860509725604\\
4	10.1136762032378\\
5	5.53755421677899\\
6	1.77438308946485\\
7	0.812932618606207\\
8	0.438976340729508\\
9	0.209183739183115\\
10	0.484255689669338\\
11	0.294677250434024\\
12	0.156485268886335\\
13	0.118988094317577\\
14	0.0894783226898153\\
15	0.0176725221191607\\
16	0.00632244326642448\\
17	0.000801363748091472\\
};
\addlegendentry{$\epsilon_{\text{true}}$}

\end{axis}
\end{tikzpicture}%
} 
	\subfloat[]{\label{fig:AVA_rand_eff}
%
%
\definecolor{mycolor1}{rgb}{0.00000,0.44700,0.74100}%
\begin{tikzpicture}

\begin{axis}[%
width=\fwidth,
height=\fheight,
at={(0\fwidth,0\fheight)},
scale only axis,
xmin=1,
xmax=17,
xlabel style={font=\color{white!15!black}},
xlabel={Iterations},
ymin=0.4,
ymax=1.6,
ylabel style={font=\color{white!15!black}},
ylabel={\texttt{eff}},
axis background/.style={fill=white},
xmajorgrids,
ymajorgrids
]
\addplot [color=white!20!black, thick]
  table[row sep=crcr]{%
1	0.709646504419979\\
2	1.05747724050268\\
3	0.997012470146273\\
4	0.940115722741877\\
5	1.50244760369346\\
6	0.810435044750712\\
7	1.42166776564107\\
8	0.778066182567579\\
9	0.824451823219507\\
10	0.910185110732025\\
11	0.990204284110803\\
12	0.776825871826319\\
13	0.458946927805223\\
14	0.634024798297823\\
15	0.920861291006627\\
16	0.965373574627563\\
17	1.24024105394485\\
};
\end{axis}
\end{tikzpicture}%
}
	\caption{Antipodal Vivaldi antenna: results for Test 3. (a) Convergence of the greedy algorithm. (b) Effectivity (\texttt{eff}).}
	\label{fig:AVA_Test3}
\end{figure}
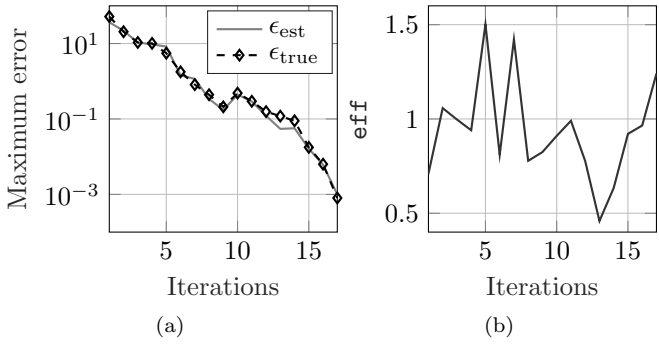
\begin{figure}[tbp]
	\centering
	\subfloat[]{\label{fig:AVA_prop_conv}
%
%
\definecolor{mycolor1}{rgb}{0.00000,0.44700,0.74100}%
\begin{tikzpicture}

\begin{axis}[%
width=\fwidth,
height=\fheight,
at={(0\fwidth,0\fheight)},
scale only axis,
xmin=1,
xmax=18,
xlabel style={font=\color{white!15!black}},
xlabel={Iterations},
ymode=log,
ymin=0.0001,
ymax=100,
yminorticks=true,
ytick = {10,0.1,0.001},
ylabel style={font=\color{white!15!black}},
ylabel={Maximum error},
axis background/.style={fill=white},
xmajorgrids,
ymajorgrids,
yminorgrids,
legend style={legend cell align=left, align=left, draw=white!15!black}
]
\addplot [color=gray,thick]
  table[row sep=crcr]{%
1	6.05412047168785\\
2	14.2830830937628\\
3	17.3764859902623\\
4	33.2098598022976\\
5	17.4135774341606\\
6	16.7398976367202\\
7	1.95142531212539\\
8	1.17044339781223\\
9	0.335753830134158\\
10	0.160115143574833\\
11	0.141614468803782\\
12	0.119378485446579\\
13	0.054334054245206\\
14	0.0330578846961133\\
15	0.0191176581116633\\
16	0.00556888736423313\\
17	0.00189789062241554\\
18	0.000836272758138832\\
};
\addlegendentry{$\epsilon_{\text{est}}$}

\addplot [color=black, dashed, mark=diamond, mark options={solid, black},thick]
  table[row sep=crcr]{%
1	51.8466903471054\\
2	51.7156165641928\\
3	34.8620491100827\\
4	34.6930985133822\\
5	17.4372875989223\\
6	16.7354991639088\\
7	1.958795776349\\
8	1.16665153028993\\
9	0.33575383013408\\
10	0.160115143574186\\
11	0.141614468800987\\
12	0.119378358264416\\
13	0.0543348643318252\\
14	0.0330578586860456\\
15	0.01911765811159\\
16	0.00556888736412149\\
17	0.00189789052388169\\
18	0.000836272758059398\\
};
\addlegendentry{$\epsilon_{\text{true}}$}

\end{axis}
\end{tikzpicture}%
} 
	\subfloat[]{\label{fig:AVA_prop_eff}
%
%
\definecolor{mycolor1}{rgb}{0.00000,0.44700,0.74100}%
\begin{tikzpicture}

\begin{axis}[%
width=\fwidth,
height=\fheight,
at={(0\fwidth,0\fheight)},
scale only axis,
xmin=1,
xmax=18,
xlabel style={font=\color{white!15!black}},
xlabel={Iterations},
ymin=0.1,
ymax=1.1,
ylabel style={font=\color{white!15!black}},
ylabel={\texttt{eff}},
axis background/.style={fill=white},
xmajorgrids,
ymajorgrids
]
\addplot [color=white!20!black, thick]
  table[row sep=crcr]{%
1	0.116769661306372\\
2	0.276185106988597\\
3	0.49843558923898\\
4	0.957246865381237\\
5	0.998640260727069\\
6	1.0002628229232\\
7	0.996237247234959\\
8	1.00325021433037\\
9	1.00000000000023\\
10	1.00000000000404\\
11	1.00000000001973\\
12	1.00000106537035\\
13	0.999985090850429\\
14	1.00000078680437\\
15	1.00000000000383\\
16	1.00000000002005\\
17	1.00000005191756\\
18	1.00000000009499\\
};
\end{axis}
\end{tikzpicture}%
}
	\caption{Antipodal Vivaldi antenna: results for Test 4. (a) Convergence of Algorithm~\ref{alg:rom_errest}. (b) Effectivity (\texttt{eff}). }
	\label{fig:AVA_Test4}
\end{figure}
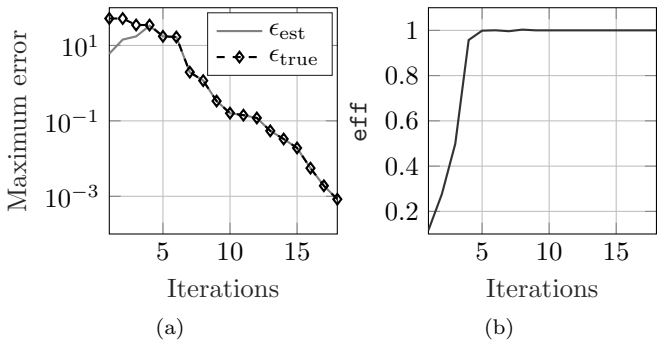
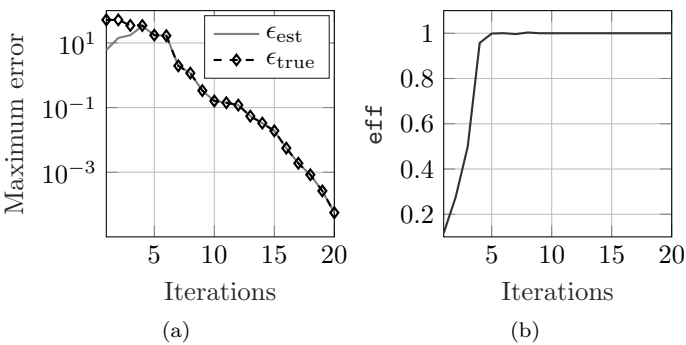
\begin{figure}[tbp]
	\centering
	\subfloat[]{\label{fig:AVA_prop_conv_smalltol}
%
%
\definecolor{mycolor1}{rgb}{0.00000,0.44700,0.74100}%
\begin{tikzpicture}

\begin{axis}[%
width=\fwidth,
height=\fheight,
at={(0\fwidth,0\fheight)},
scale only axis,
xmin=1,
xmax=20,
xlabel style={font=\color{white!15!black}},
xlabel={Iterations},
ymode=log,
ymin=1e-05,
ymax=100,
yminorticks=true,
ytick = {10,0.1,0.001},
ylabel style={font=\color{white!15!black}},
ylabel={Maximum error},
axis background/.style={fill=white},
xmajorgrids,
ymajorgrids,
yminorgrids,
legend style={legend cell align=left, align=left, draw=white!15!black}
]
\addplot [color=gray, thick]
  table[row sep=crcr]{%
1	6.05412047168785\\
2	14.2830830937628\\
3	17.3764859902623\\
4	33.2098598022976\\
5	17.4135774341606\\
6	16.7398976367202\\
7	1.95142531212539\\
8	1.17044339781223\\
9	0.335753830134158\\
10	0.160115143574833\\
11	0.141614468803782\\
12	0.119378485446579\\
13	0.054334054245206\\
14	0.0330578846961133\\
15	0.0191176581116633\\
16	0.00556888736423313\\
17	0.00189789062241554\\
18	0.000836272758138832\\
19	0.000264305370361026\\
20	5.6045231542209e-05\\
};
\addlegendentry{$\epsilon_{\text{est}}$}

\addplot [color=black, dashed, mark=diamond, mark options={solid, black},thick]
  table[row sep=crcr]{%
1	51.8466903471054\\
2	51.7156165641928\\
3	34.8620491100827\\
4	34.6930985133822\\
5	17.4372875989223\\
6	16.7354991639088\\
7	1.958795776349\\
8	1.16665153028993\\
9	0.33575383013408\\
10	0.160115143574186\\
11	0.141614468800987\\
12	0.119378358264416\\
13	0.0543348643318252\\
14	0.0330578586860456\\
15	0.01911765811159\\
16	0.00556888736412149\\
17	0.00189789052388169\\
18	0.000836272758059398\\
19	0.00026430538021215\\
20	5.60452315984715e-05\\
};
\addlegendentry{$\epsilon_{\text{true}}$}

\end{axis}
\end{tikzpicture}%
} 
	\subfloat[]{\label{fig:AVA_prop_eff_smalltol}
%
%
\definecolor{mycolor1}{rgb}{0.00000,0.44700,0.74100}%
\begin{tikzpicture}

\begin{axis}[%
width=\fwidth,
height=\fheight,
at={(0\fwidth,0\fheight)},
scale only axis,
xmin=1,
xmax=20,
xlabel style={font=\color{white!15!black}},
xlabel={Iterations},
ymin=0.1,
ymax=1.1,
ylabel style={font=\color{white!15!black}},
ylabel={\texttt{eff}},
axis background/.style={fill=white},
xmajorgrids,
ymajorgrids
]
\addplot [color=white!20!black, thick]
  table[row sep=crcr]{%
1	0.116769661306372\\
2	0.276185106988597\\
3	0.49843558923898\\
4	0.957246865381237\\
5	0.998640260727069\\
6	1.0002628229232\\
7	0.996237247234959\\
8	1.00325021433037\\
9	1.00000000000023\\
10	1.00000000000404\\
11	1.00000000001973\\
12	1.00000106537035\\
13	0.999985090850429\\
14	1.00000078680437\\
15	1.00000000000383\\
16	1.00000000002005\\
17	1.00000005191756\\
18	1.00000000009499\\
19	0.999999962728251\\
20	0.999999998996122\\
};
\end{axis}
\end{tikzpicture}%
}
	\caption{Antipodal Vivaldi antenna: results for Test 4 with tolerance $\texttt{tol} = 10^{-4}$. (a) Convergence of Algorithm~\ref{alg:rom_errest}. (b) Effectivity (\texttt{eff}). }
	\label{fig:AVA_Test4_smalltol}
\end{figure}
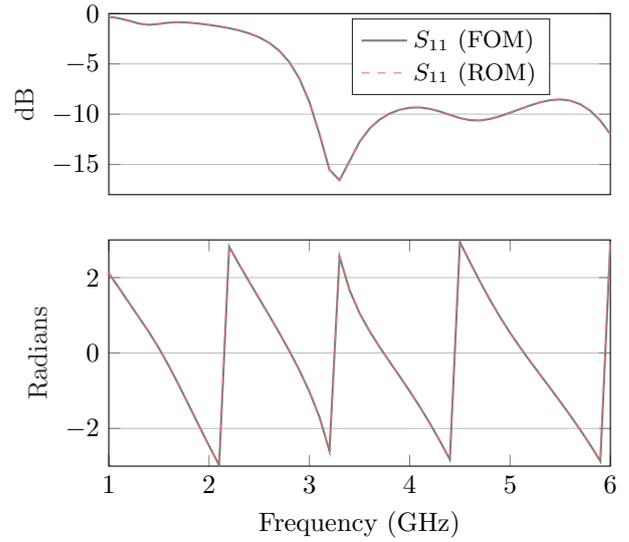
\begin{figure}[t!]
	\centering
%
%
\begin{tikzpicture}

\begin{axis}[%
width=2.2\fwidth,
height=0.8\fheight,
at={(0\fwidth,0.8\fheight)},
scale only axis,
xmin=1,
xmax=6,
xtick={\empty},
ymin=-18,
ymax=0,
ylabel style={font=\color{white!15!black}},
ylabel={dB},
axis background/.style={fill=white},
xmajorgrids,
ymajorgrids,
legend pos=south east,
legend style={at={(0.88, 0.55)},legend cell align=left, align=left, draw=white!15!black, font=\small}
]
\addplot [color=gray, thick]
table[row sep=crcr]{%
1	-0.30357404889455\\
1.1	-0.461711773739517\\
1.2	-0.703145302536909\\
1.3	-0.986902799637892\\
1.4	-1.11741821282661\\
1.5	-1.02000037267867\\
1.6	-0.895523385806415\\
1.7	-0.855649146276355\\
1.8	-0.895991918103159\\
1.9	-0.989638762008902\\
2	-1.11643027964236\\
2.1	-1.26818672160964\\
2.2	-1.44874072657343\\
2.3	-1.6724611158722\\
2.4	-1.96361148998657\\
2.5	-2.35820798842476\\
2.6	-2.90899888683366\\
2.7	-3.69352833348248\\
2.8	-4.82434400037417\\
2.9	-6.45832771796033\\
3	-8.79292677208177\\
3.1	-11.9885785880011\\
3.2	-15.5368539586489\\
3.3	-16.5584690699273\\
3.4	-14.6618126723269\\
3.5	-12.7489600953815\\
3.6	-11.4204325231568\\
3.7	-10.5312913081605\\
3.8	-9.9350381814926\\
3.9	-9.55456565802618\\
4	-9.35873222680532\\
4.1	-9.33522180309748\\
4.2	-9.46640535724765\\
4.3	-9.72918682210344\\
4.4	-10.0642430484377\\
4.5	-10.3854140085405\\
4.6	-10.5948741026839\\
4.7	-10.6249240040576\\
4.8	-10.4739595609277\\
4.9	-10.1967937832843\\
5	-9.85581448136167\\
5.1	-9.49164178846555\\
5.2	-9.13439981756681\\
5.3	-8.82226822234662\\
5.4	-8.60395151207669\\
5.5	-8.53074867079819\\
5.6	-8.65028937916267\\
5.7	-9.00893125204908\\
5.8	-9.65406324706281\\
5.9	-10.6344562831643\\
6	-11.989410366968\\
};
\addlegendentry{$S_{11}$ (FOM)}

\addplot [color=red!60, dashed]
table[row sep=crcr]{%
1	-0.303574048894484\\
1.1	-0.461711774191786\\
1.2	-0.703145302537291\\
1.3	-0.986902803936598\\
1.4	-1.11741821486401\\
1.5	-1.02000037267846\\
1.6	-0.895523385805643\\
1.7	-0.855649128184063\\
1.8	-0.895991753908826\\
1.9	-0.989638257906632\\
2	-1.11642947082706\\
2.1	-1.2681860937019\\
2.2	-1.44874102915842\\
2.3	-1.67246283201831\\
2.4	-1.96361436159163\\
2.5	-2.35821101385657\\
2.6	-2.90900093013666\\
2.7	-3.69352899879204\\
2.8	-4.82434400037382\\
2.9	-6.45832771796011\\
3	-8.79292679555424\\
3.1	-11.9885781683809\\
3.2	-15.5368535855661\\
3.3	-16.558469069933\\
3.4	-14.6618123329301\\
3.5	-12.7489595881543\\
3.6	-11.4204324354811\\
3.7	-10.5312913081615\\
3.8	-9.93503826946637\\
3.9	-9.55456571482995\\
4	-9.35873222680479\\
4.1	-9.33522154848552\\
4.2	-9.46640535724702\\
4.3	-9.72918671884576\\
4.4	-10.0642430545778\\
4.5	-10.38541400854\\
4.6	-10.5948741281354\\
4.7	-10.6249240073982\\
4.8	-10.473959560928\\
4.9	-10.1967937752831\\
5	-9.8558144801097\\
5.1	-9.49164178846581\\
5.2	-9.13439981843296\\
5.3	-8.82226822234705\\
5.4	-8.60395151184846\\
5.5	-8.53074867079879\\
5.6	-8.65028937916379\\
5.7	-9.00893123028975\\
5.8	-9.65406324706213\\
5.9	-10.6344562755639\\
6	-11.9894103669685\\
};
\addlegendentry{$S_{11}$ (ROM)}
\end{axis}

\begin{axis}[%
width=2.2\fwidth,
height=\fheight,
at={(0\fwidth,-0.4\fheight)},
scale only axis,
xmin=1,
xmax=6,
xlabel={Frequency (GHz)},
ymin=-3,
ymax=3,
ylabel style={font=\color{white!15!black}},
ylabel={Radians},
axis background/.style={fill=white},
ymajorgrids,
legend style={legend cell align=left, align=left, draw=white!15!black}
]
\addplot [color=gray, thick]
table[row sep=crcr]{%
1	2.1153631297756\\
1.1	1.72968357447136\\
1.2	1.33868125412734\\
1.3	0.948736544125576\\
1.4	0.556096103027613\\
1.5	0.12855724843081\\
1.6	-0.348165920886171\\
1.7	-0.860840732362728\\
1.8	-1.39334052262158\\
1.9	-1.93097694669521\\
2	-2.46078955908259\\
2.1	-2.97358051259805\\
2.2	2.8178076667996\\
2.3	2.34597897479597\\
2.4	1.88937627731321\\
2.5	1.44089301883148\\
2.6	0.991990386571225\\
2.7	0.533157055317656\\
2.8	0.0538629693577754\\
2.9	-0.458709278593247\\
3	-1.02579793116384\\
3.1	-1.69812755535679\\
3.2	-2.59871366626363\\
3.3	2.5591430172667\\
3.4	1.66826157460801\\
3.5	1.05504981192508\\
3.6	0.575650963376899\\
3.7	0.158255580867276\\
3.8	-0.232806816104939\\
3.9	-0.618702598225071\\
4	-1.01315440783013\\
4.1	-1.42560318893016\\
4.2	-1.86248351207524\\
4.3	-2.32881636396779\\
4.4	-2.82428029247335\\
4.5	2.94049288669038\\
4.6	2.4135386398068\\
4.7	1.89723218055437\\
4.8	1.40789605335139\\
4.9	0.95397087600132\\
5	0.536084515279653\\
5.1	0.148570836972425\\
5.2	-0.218233962556002\\
5.3	-0.574726446807041\\
5.4	-0.929480086100745\\
5.5	-1.28863279346529\\
5.6	-1.65671381095725\\
5.7	-2.03847520895843\\
5.8	-2.43926662041283\\
5.9	-2.86777609794454\\
6	2.94484047152479\\
};

\addplot [color=red!60, dashed]
table[row sep=crcr]{%
1	2.11536312977549\\
1.1	1.72968357416378\\
1.2	1.33868125412742\\
1.3	0.948736544122681\\
1.4	0.556096103351939\\
1.5	0.128557248430795\\
1.6	-0.348165920886053\\
1.7	-0.860840732186808\\
1.8	-1.39334052932515\\
1.9	-1.93097699588811\\
2	-2.46078971434646\\
2.1	-2.97358082235082\\
2.2	2.81780722259919\\
2.3	2.34597849753349\\
2.4	1.8893759019107\\
2.5	1.44089283052446\\
2.6	0.991990360825783\\
2.7	0.53315708149215\\
2.8	0.0538629693575972\\
2.9	-0.458709278593273\\
3	-1.02579786770111\\
3.1	-1.69812744916797\\
3.2	-2.59871361491145\\
3.3	2.55914301726784\\
3.4	1.66826159284882\\
3.5	1.05504979763634\\
3.6	0.575650943995943\\
3.7	0.158255580867152\\
3.8	-0.232806817512998\\
3.9	-0.618702590958502\\
4	-1.01315440782999\\
4.1	-1.42560318627037\\
4.2	-1.86248351207507\\
4.3	-2.32881637446135\\
4.4	-2.8242802995272\\
4.5	2.94049288669035\\
4.6	2.41353864104235\\
4.7	1.89723218280222\\
4.8	1.40789605335125\\
4.9	0.953970875697634\\
5	0.53608451473762\\
5.1	0.148570836972413\\
5.2	-0.218233962553151\\
5.3	-0.574726446807097\\
5.4	-0.929480085987471\\
5.5	-1.28863279346541\\
5.6	-1.65671381095735\\
5.7	-2.03847520850751\\
5.8	-2.43926662041288\\
5.9	-2.86777610922036\\
6	2.94484047152485\\
};
\end{axis}

\end{tikzpicture}%
	\caption{Antipodal Vivaldi antenna: scattering parameter responses of FOM and ROM computed from Test 4. Top: The magnitude of the scattering parameter response $|S_{11}|$; Bottom: The phase $\angle(S_{11})$.}
	\label{fig:AVA_Sparameters}
\end{figure}

\subsection{Example IV : Dielectric Resonator Antenna}
The next example we consider is the model of a dielectric resonator antenna \cite{delZ07} shown in Fig.~\ref{fig:DRA_Geometry}. Reducing the metallization part, such antennas tend to have low losses and are light weight \cite{LiL05}. Among the six examples considered, this model has the largest dimension with $n = 484,294$. The frequency band of interest spans $[2.5, \, 4.5]$ GHz. For the training set, we sample 41 points uniformly. The tolerance for the greedy algorithms is $\texttt{tol} = 10^{-3}$.%
\subsubsection{Test 2. Residual norm as a heuristic error estimator}
The residual error estimator serves as a fairly good surrogate for the true error, as can be seen in Fig.~\ref{fig:DRA_Test2}. The algorithm converges in $12$ iterations to a ROM of size $r = 24$. The runtime for the algorithm until convergence is around $49$ minutes.\\
\subsubsection{Test 3. Randomized error estimator}
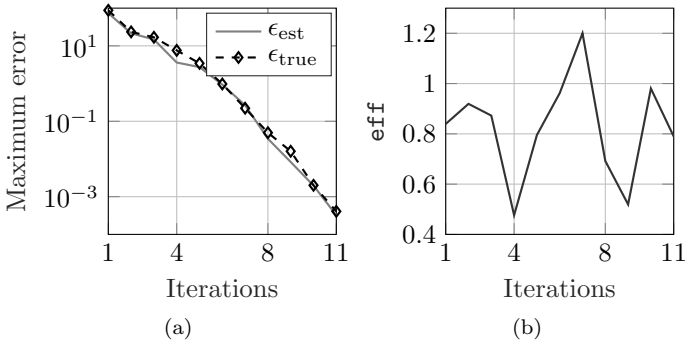
\begin{figure}[t!]
	\centering
	\subfloat[]{\label{fig:DRA_rand_conv}
%
%
\definecolor{mycolor1}{rgb}{0.00000,0.44700,0.74100}%
\begin{tikzpicture}

\begin{axis}[%
width=\fwidth,
height=\fheight,
at={(0\fwidth,0\fheight)},
scale only axis,
xmin=1,
xmax=11,
xtick = {1,4,8,11},
xlabel style={font=\color{white!15!black}},
xlabel={Iterations},
ymode=log,
ymin=0.0001,
ymax=100,
yminorticks=true,
ytick = {10,0.1,0.001},
ylabel style={font=\color{white!15!black}},
ylabel={Maximum error},
axis background/.style={fill=white},
xmajorgrids,
ymajorgrids,
yminorgrids,
legend style={legend cell align=left, align=left, draw=white!15!black}
]
\addplot [color=gray,thick]
  table[row sep=crcr]{%
1	73.5641213855276\\
2	21.564679488665\\
3	14.7181224074449\\
4	3.61282534227511\\
5	2.73437326876592\\
6	0.937403991415032\\
7	0.266259386745253\\
8	0.034431491618664\\
9	0.00832115684256185\\
10	0.0019661608938589\\
11	0.000322546503832051\\
};
\addlegendentry{$\epsilon_{\text{est}}$}

\addplot [color=black, dashed, mark=diamond, mark options={solid, black},thick]
  table[row sep=crcr]{%
1	87.6391010366013\\
2	23.4574315030513\\
3	16.8674362413195\\
4	7.5851231833021\\
5	3.43766455110699\\
6	0.974579291153039\\
7	0.221892734288256\\
8	0.0497439863620597\\
9	0.0160253396962215\\
10	0.00200553972839836\\
11	0.000408511322166446\\
};
\addlegendentry{$\epsilon_{\text{true}}$}

\end{axis}
\end{tikzpicture}%
} 
	\subfloat[]{\label{fig:DRA_rand_eff}
%
%
\definecolor{mycolor1}{rgb}{0.00000,0.44700,0.74100}%
\begin{tikzpicture}

\begin{axis}[%
width=\fwidth,
height=\fheight,
at={(0\fwidth,0\fheight)},
scale only axis,
xmin=1,
xmax=11,
xtick = {1,4,8,11},
xlabel style={font=\color{white!15!black}},
xlabel={Iterations},
ymin=0.4,
ymax=1.3,
ylabel style={font=\color{white!15!black}},
ylabel={\texttt{eff}},
axis background/.style={fill=white},
xmajorgrids,
ymajorgrids
]
\addplot [color=white!20!black, thick]
  table[row sep=crcr]{%
1	0.839398402258879\\
2	0.919311199346778\\
3	0.872576140017677\\
4	0.47630410936877\\
5	0.795415965727489\\
6	0.961855028035713\\
7	1.19994639571822\\
8	0.692173951803053\\
9	0.519249950409713\\
10	0.980364969099413\\
11	0.789565640730591\\
};
\end{axis}
\end{tikzpicture}%
}
	\caption{Dielectric resonator antenna: results for Test 3. (a) Convergence of the greedy algorithm. (b) Effectivity (\texttt{eff}).}
	\label{fig:DRA_Test3}
\end{figure}

We set $K = 5$, $\texttt{tol}_{rd} = 1$ to determine $V_{rd}$. Once again, the random number generator \texttt{MersenneTwister} was used with the seed set to $1$ to draw the random vectors. As for the previous example, this process takes up a large time, requiring almost $113$ minutes. Fig.~\ref{fig:DRA_rand_conv} shows the convergence of the greedy algorithm while Fig.~\ref{fig:DRA_rand_eff} displays the effectivity of the estimator. Evidently, the randomized error estimator shows a mean effectivity around one, but with a tendency to slightly underestimate the true error. The resulting ROM of size $r = 22$ is achieved in $11$ iterations.\\
\subsubsection{Test 4. Proposed error estimator}
For the final test, we apply Algorithm 1 to reduce the model of the dielectric resonator antenna. Fig.~\ref{fig:DRA_Test4} shows the procedure results in a ROM with size $r = 22$, as in Test 3. However, unlike the randomized error estimator, there is no underestimation of the true error. In fact, the effectivity is almost unity, for most of the $11$ iterations. Also, the time required for Algorithm 1 is about $90$ minutes. As done for the previous examples, we also plot the scattering parameter response for this antenna in Fig.~\ref{fig:DRA_Sparameters}. We see that the response obtained from the ROM is almost exactly matching the results obtained using the full order model, thus showing the good performance of the ROM.
\begin{figure}[tbp]
	\centering
	\includegraphics[width=\linewidth]{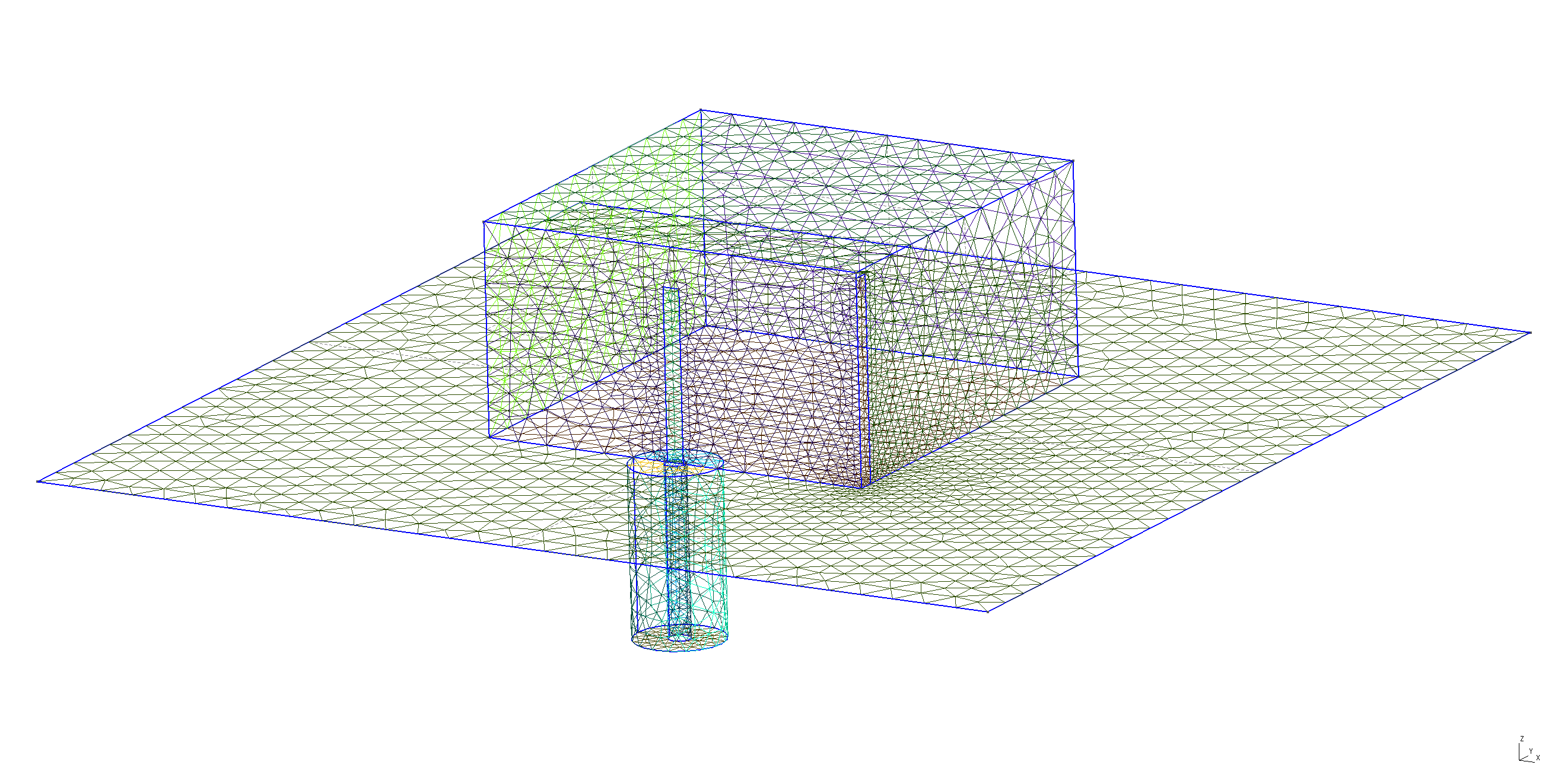}
	\caption{Dielectric resonator antenna designed in \cite{LiL05}.}
	\label{fig:DRA_Geometry}
\end{figure}
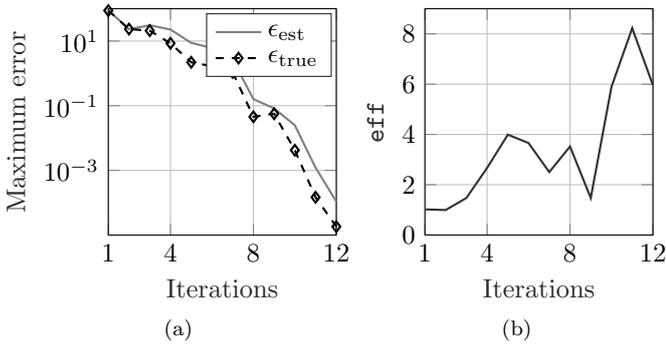
\begin{figure}[tbp]
	\centering
	\subfloat[]{\label{fig:DRA_resd_conv}
%
%
\definecolor{mycolor1}{rgb}{0.00000,0.44700,0.74100}%
\begin{tikzpicture}

\begin{axis}[%
width=\fwidth,
height=\fheight,
at={(0\fwidth,0\fheight)},
scale only axis,
xmin=1,
xmax=12,
xtick = {1,4,8,12},
xlabel style={font=\color{white!15!black}},
xlabel={Iterations},
ymode=log,
ymin=1e-05,
ymax=100,
yminorticks=true,
ytick = {10,0.1,0.001},
ylabel style={font=\color{white!15!black}},
ylabel={Maximum error},
axis background/.style={fill=white},
xmajorgrids,
ymajorgrids,
yminorgrids,
legend style={legend cell align=left, align=left, draw=white!15!black}
]
\addplot [color=gray,thick]
  table[row sep=crcr]{%
1	89.4290727368796\\
2	23.2834893436026\\
3	30.8804905815324\\
4	22.7502038634471\\
5	8.78367187208118\\
6	6.10259678486563\\
7	2.67720856277902\\
8	0.15982572766022\\
9	0.0828648555057897\\
10	0.0247291412385299\\
11	0.00121007119305807\\
12	0.000108170513108896\\
};
\addlegendentry{$\epsilon_{\text{est}}$}

\addplot [color=black, dashed, mark=diamond, mark options={solid, black},thick]
  table[row sep=crcr]{%
1	87.6391010366013\\
2	23.4574315030513\\
3	20.974019751632\\
4	8.51902215287302\\
5	2.20073229003829\\
6	1.6681485388361\\
7	1.06998547482376\\
8	0.0453954949931806\\
9	0.0562745368870652\\
10	0.00419444164375661\\
11	0.000147093375344423\\
12	1.81510331709896e-05\\
};
\addlegendentry{$\epsilon_{\text{true}}$}

\end{axis}
\end{tikzpicture}%
} 
	\subfloat[]{\label{fig:DRA_resd_eff}
%
%
\definecolor{mycolor1}{rgb}{0.00000,0.44700,0.74100}%
\begin{tikzpicture}

\begin{axis}[%
width=\fwidth,
height=\fheight,
at={(0\fwidth,0\fheight)},
scale only axis,
xmin=1,
xmax=12,
xtick = {1,4,8,12},
xlabel style={font=\color{white!15!black}},
xlabel={Iterations},
ymin=0,
ymax=9,
ylabel style={font=\color{white!15!black}},
ylabel={\texttt{eff}},
axis background/.style={fill=white},
xmajorgrids,
ymajorgrids
]
\addplot [color=white!20!black, thick]
  table[row sep=crcr]{%
1	1.02042435030833\\
2	0.992584773851902\\
3	1.47232104037327\\
4	2.67051821854632\\
5	3.99124959989039\\
6	3.65830538635579\\
7	2.50209804317203\\
8	3.5207398373832\\
9	1.47251066094222\\
10	5.89569323853606\\
11	8.22655126530786\\
12	5.95946864786643\\
};
\end{axis}
\end{tikzpicture}%
}
	\caption{Dielectric resonator antenna: results for Test 2. (a) Convergence of the greedy algorithm. (b) Effectivity (\texttt{eff}).}
	\label{fig:DRA_Test2}
\end{figure}

\begin{figure}[tbp]
	\centering
	\subfloat[]{\label{fig:DRA_prop_conv}
%
%
\definecolor{mycolor1}{rgb}{0.00000,0.44700,0.74100}%
\begin{tikzpicture}

\begin{axis}[%
width=\fwidth,
height=\fheight,
at={(0\fwidth,0\fheight)},
scale only axis,
xmin=1,
xmax=11,
xtick = {1,4,8,11},
xlabel style={font=\color{white!15!black}},
xlabel={Iterations},
ymode=log,
ymin=0.0001,
ymax=100,
yminorticks=true,
ytick = {10,0.1,0.001},
ylabel style={font=\color{white!15!black}},
ylabel={Maximum error},
axis background/.style={fill=white},
xmajorgrids,
ymajorgrids,
yminorgrids,
legend style={legend cell align=left, align=left, draw=white!15!black}
]
\addplot [color=gray,thick]
  table[row sep=crcr]{%
1	31.0837748985066\\
2	96.844939216515\\
3	19.1240208920798\\
4	4.85143267339865\\
5	2.7570891406425\\
6	1.0253683995871\\
7	0.236589764596759\\
8	0.0600668967700996\\
9	0.0191718486349041\\
10	0.00172582716156706\\
11	0.000206635626991782\\
};
\addlegendentry{$\epsilon_{\text{est}}$}

\addplot [color=black, dashed, mark=diamond, mark options={solid, black},thick]
  table[row sep=crcr]{%
1	87.6391010366013\\
2	96.3534623262932\\
3	19.3912679765907\\
4	4.85178367479144\\
5	2.75826632759402\\
6	1.0253656553795\\
7	0.23658976459676\\
8	0.0600668981867746\\
9	0.0191718483967363\\
10	0.00172582716150862\\
11	0.000206635628390866\\
};
\addlegendentry{$\epsilon_{\text{true}}$}

\end{axis}
\end{tikzpicture}%
} 
	\subfloat[]{\label{fig:DRA_prop_eff}
%
%
\definecolor{mycolor1}{rgb}{0.00000,0.44700,0.74100}%
\begin{tikzpicture}

\begin{axis}[%
width=\fwidth,
height=\fheight,
at={(0\fwidth,0\fheight)},
scale only axis,
xmin=1,
xmax=11,
xtick = {1,4,8,11},
xlabel style={font=\color{white!15!black}},
xlabel={Iterations},
ymin=0.8,
ymax=1.05,
ylabel style={font=\color{white!15!black}},
ylabel={\texttt{eff}},
axis background/.style={fill=white},
xmajorgrids,
ymajorgrids
]
\addplot [color=white!20!black, thick]
  table[row sep=crcr]{%
1	0.989744943613992\\
2	0.805582889551219\\
3	0.9933239131711\\
4	0.999850728894765\\
5	0.999693331733628\\
6	1.00000154099014\\
7	0.999999999999245\\
8	0.999999985751543\\
9	0.999999999657125\\
10	1.00000000005687\\
11	0.999999989239612\\
};
\end{axis}
\end{tikzpicture}%
}
	\caption{Dielectric resonator antenna: results for Test 4. (a) Convergence of Algorithm~\ref{alg:rom_errest}. (b) Effectivity (\texttt{eff}).}
	\label{fig:DRA_Test4}
\end{figure}
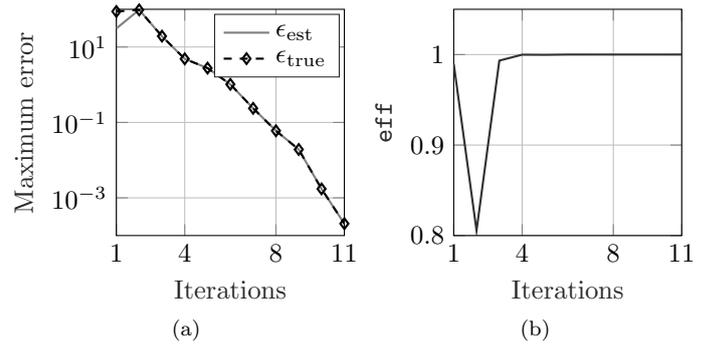
\begin{figure}[t!]
	\centering
	\input{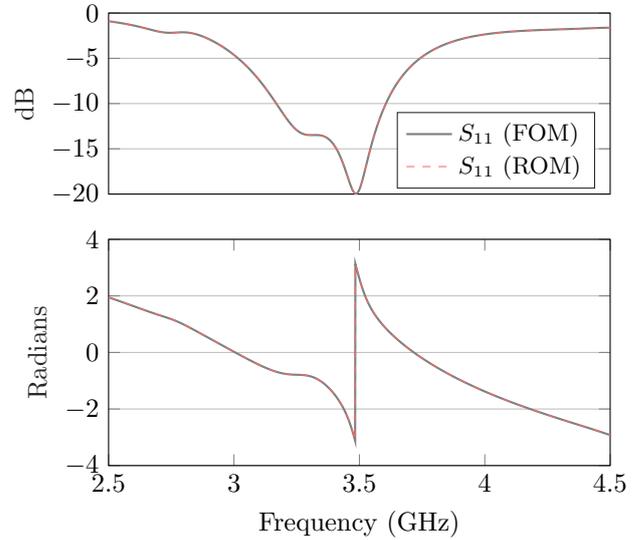}
	\caption{Dielectric resonator antenna: scattering parameter responses of FOM and ROM computed from Test 4. Top: The magnitude of the scattering parameter response $|S_{11}|$; Bottom: The phase $\angle(S_{11})$.}
	\label{fig:DRA_Sparameters}
\end{figure}

\subsection{Example V : Coax-fed Dielectric Resonator Filter}
\begin{figure}[tbp]
	\centering
	\includegraphics[width=0.85\linewidth]{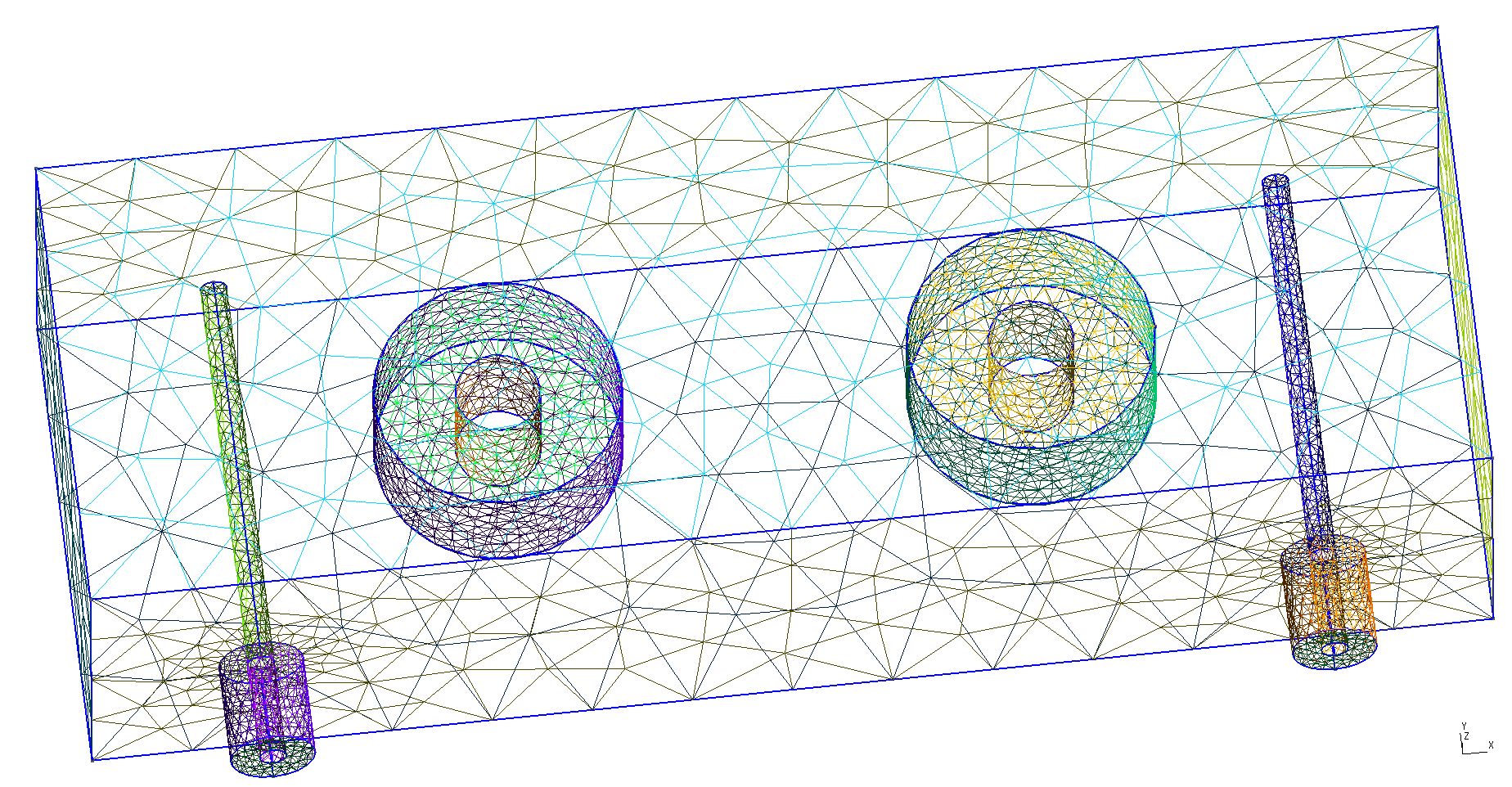} 
	\caption{Geometry of the coax-fed dielectric resonator filter from \cite{BraL97}.}
	\label{fig:Brauer_Geometry}
\end{figure}

Next, we consider the coax-fed dielectric resonator filter originally proposed in \cite{BraL97}. The geometry of the filter is shown in Fig.~\ref{fig:Brauer_Geometry}. The filter consists of two cylindrical dielectric resonators, each having a concentric hole. The discretized model is of dimension $n = 154,066$. We consider a wide frequency range of interest, i.e., $f \in [4.0, 12.0]$ GHz. This example is quite challenging owing to the presence of multiple resonances in its frequency range. The tolerance for the ROM is set as $\texttt{tol} = 10^{-4}$ and the training set consists of $2000$ uniformly-sampled parameters from the parameter range.

Due to the large number of samples in the training set, we refrain from showing the true error convergence for this example. This is due to the significant computational cost and storage requirements associated with determining the true solutions. Note that the true error is computed only for comparison, and its computation is always avoided in practical applications anyway.
%
%
\subsubsection{Test 2. Residual norm as a heuristic error estimator}
By using the norm of the residual $\| r(s) \|$ as the error estimator in Algorithm \ref{alg:rom_errest}, we obtain a ROM of dimension $r = 127$. The convergence of the greedy algorithm is illustrated in Fig.~\ref{fig:Brauer_Test234_conv} (dotted line). The total number of iterations required is $64$ and the time taken to converge is $12$ minutes and $32$ seconds.\\


\subsubsection{Test 3. Randomized error estimator}
Applying the randomized error estimator from \cite{morSemZP18} leads to a ROM with slightly smaller dimension of $r = 124$, with the greedy algorithm taking $14$ minutes and $62$ iterations to converge as seen in Fig.~\ref{fig:Brauer_Test234_conv} (dashed line). A large part of the offline time is spent in constructing the matrix $V_{rd}$, which took about $204$ minutes. To determine $V_{rd}$, we set $K=10$ and $\texttt{tol}_{rd} = 10$. The $K$ random vectors were drawn using $\texttt{mvnrnd}$ in MATLAB\textsuperscript \textregistered$\,$ using the $\texttt{simdTwister}$ random number generator, with the seed set to $10$.\\


\subsubsection{Test 4. Proposed error estimator}
The results of using the proposed error estimator in Algorithm \ref{alg:rom_errest} is shown in Fig.~\ref{fig:Brauer_Test234_conv}. The greedy algorithm needs $61$ iterations to achieve the desired tolerance of $10^{-4}$, taking around $30$ minutes to converge. The resulting ROM has dimension $r = 122$. The proposed approach thus results in a ROM having the smallest dimension among the three tests. 

While it indeed takes a longer time to converge when compared to Test 2, the proposed approach enjoys the benefits of having a guaranteed accuracy of the estimated error and a ROM of smaller dimension.

\begin{figure}[t!]
	\centering
	\includegraphics[scale=1]{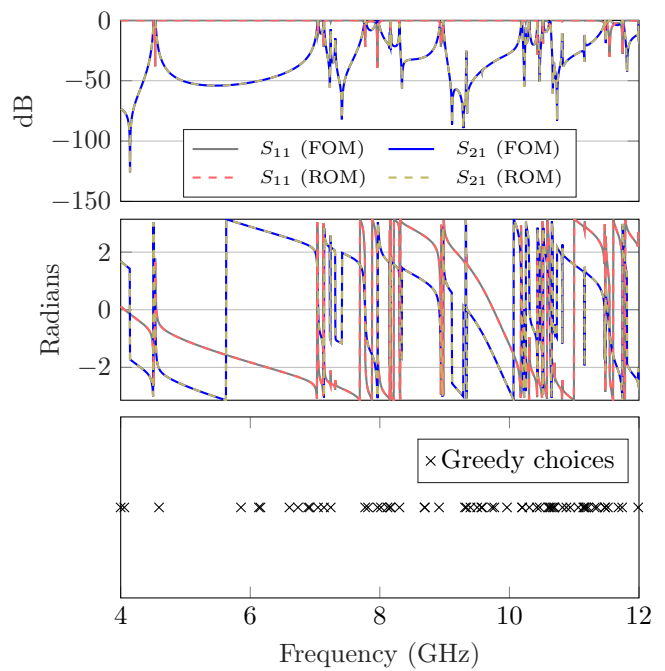}
	\caption{Coax-fed dielectric resonator filter: scattering parameter responses of FOM and ROM computed from Test 4. Top: The magnitudes of the scattering parameter responses $|S_{11}|$, $|S_{21}|$; Middle: The phases $\angle(S_{11}), \angle(S_{21})$; Bottom: The greedy parameters $\mu^{*}$ selected by Algorithm~\ref{alg:rom_errest} using the proposed error estimator.}
	\label{fig:Brauer_Test4_gp}
\end{figure}

The scattering parameter computed from the FOM and the ROM obtained from Test 4 are compared in Fig.~\ref{fig:Brauer_Test4_gp}. The magnitude and the phase are plotted separately. The ROM offers an excellent match with the true scattering response for this example. In the same figure, we also show the locations of the greedy parameters $\mu^{*}$ picked during Test 4. More samples are concentrated in the higher frequency regions, owing to the presence of a significant number of resonances there. It is also noteworthy to see that Algorithm~\ref{alg:rom_errest} avoids sampling in the range $[4.6, 7]$ GHz where no resonances appear.

\setlength\fheight{4cm}
\setlength\fwidth{5cm}
\begin{figure}[tbp]
	\centering
%
%
\definecolor{mycolor1}{rgb}{0.00000,0.44700,0.74100}%
\begin{tikzpicture}

\begin{axis}[%
width=\fwidth,
height=\fheight,
at={(0\fwidth,0\fheight)},
scale only axis,
xmin=1,
xmax=64,
xtick = {1,16,32,48,64},
xlabel style={font=\color{white!15!black}},
xlabel={Iterations},
ymode=log,
ymin=1e-06,
ymax=1000000,
yminorticks=true,
ytick = {1000,0.1,0.0001},
ylabel style={font=\color{white!15!black}},
ylabel={Maximum error},
axis background/.style={fill=white},
xmajorgrids,
ymajorgrids,
yminorgrids,
legend style={at={(0.4,0.05)}, anchor=south, legend cell align=center, align=left, font=\small, draw=white!15!black}
]
\addplot [color=gray, densely dotted, thick]
  table[row sep=crcr]{%
1	401.769276551631\\
2	4.05637741645276\\
3	3.06232480717181\\
4	607.029056130561\\
5	80.850398632016\\
6	623.421491555043\\
7	680.442217936453\\
8	758.884631530396\\
9	981.868763072518\\
10	2448.16031961867\\
11	533.429709308049\\
12	1727.59430862444\\
13	322.8987035511\\
14	1064.61786642987\\
15	719.999473707311\\
16	2282.8478331932\\
17	1803.80921416728\\
18	6989.10182312527\\
19	1217.74965340293\\
20	3590.19949660759\\
21	1456.43000988404\\
22	219.866147906132\\
23	90.8611903100913\\
24	280.636072664236\\
25	121.171217390892\\
26	52.1823832081257\\
27	43.8660549624065\\
28	133.731317826061\\
29	179.187246240133\\
30	328.048526828812\\
31	33.8760883724875\\
32	209.780132507419\\
33	31.4934872189341\\
34	24.2287829919088\\
35	47.7222054294195\\
36	18.3575996651583\\
37	26.3593253091796\\
38	329.420917956264\\
39	5.70486850629294\\
40	29.4684195835842\\
41	7.43512608809756\\
42	16.5895415906041\\
43	7.52458448519638\\
44	2.22245938170174\\
45	33.9880231601011\\
46	7.99565544727391\\
47	3.77555017418152\\
48	6.38240391257385\\
49	0.164994422331851\\
50	0.850861613535414\\
51	0.321037183017577\\
52	0.115552183619735\\
53	0.135916805748663\\
54	2.02142707896471\\
55	0.128832457354777\\
56	0.133690879178475\\
57	0.462058451373165\\
58	0.134651909589626\\
59	0.0114040673640097\\
60	0.00481785563294526\\
61	0.00189502118072348\\
62	0.000399321392137649\\
63	0.000549013806890364\\
64	9.48622253422035e-05\\
};
\addlegendentry{$\epsilon_{\text{est}}$ (residual)}
%
\addplot [color=gray,dashed,thick]
table[row sep=crcr]{%
	1	308.052729595347\\
	2	281.21862733545\\
	3	249.206943854356\\
	4	404.95026081581\\
	5	943.359080257565\\
	6	952.340385450657\\
	7	597.723429903324\\
	8	3919.31435681011\\
	9	277.488299543748\\
	10	255.34331147179\\
	11	174.576775266058\\
	12	130.052375306145\\
	13	165.238219595343\\
	14	55.8634595893923\\
	15	42.8291244139205\\
	16	35.7732073819863\\
	17	141.121783994176\\
	18	384.216847800938\\
	19	43.8703722400201\\
	20	30.1820677432171\\
	21	77.9640165766626\\
	22	491.348026765685\\
	23	89.1790981153603\\
	24	169.248266592425\\
	25	28.1438469241952\\
	26	6.2665728160689\\
	27	4.19496657948585\\
	28	3.77524102169824\\
	29	6.30263080004265\\
	30	2.34676067416671\\
	31	2.98571895218128\\
	32	1.44611783996298\\
	33	2.32176512991901\\
	34	17.2142257132486\\
	35	21.4363896755956\\
	36	35.7631815472485\\
	37	1.78484650457596\\
	38	0.790051562134285\\
	39	0.438786573154738\\
	40	2.70074453650692\\
	41	4.62166753159771\\
	42	0.161465336868609\\
	43	0.348301774321629\\
	44	0.345335508511465\\
	45	0.117128776061192\\
	46	0.329620205042033\\
	47	0.098181047368398\\
	48	1.35341930775851\\
	49	0.608291150530118\\
	50	1.67551550703925\\
	51	7.25791227959752\\
	52	0.058585918665325\\
	53	0.114866715964485\\
	54	0.0275013541290725\\
	55	0.0102980323831071\\
	56	0.00853431272706023\\
	57	0.018412284601825\\
	58	0.136325108393228\\
	59	0.0301645612671093\\
	60	0.000931397397117141\\
	61	0.000163787570144542\\
	62	1.7283179155693e-05\\
};
\addlegendentry{$\epsilon_{\text{est}}$ (randomized)}
%
\addplot [color=black, thick]
table[row sep=crcr]{%
	1	223.130893253267\\
	2	24.541390817623\\
	3	254.096415678352\\
	4	369.00737710154\\
	5	37133.5693189637\\
	6	1023.02865140243\\
	7	416.098748966498\\
	8	3343.17033528364\\
	9	141278.799050034\\
	10	96.5536935116767\\
	11	166.542470296353\\
	12	419.155304479279\\
	13	4463.61473952867\\
	14	8608.66702381543\\
	15	120.398036954114\\
	16	2524.08553852473\\
	17	735.202912726492\\
	18	159.308468303849\\
	19	104.985218366732\\
	20	38.399070456855\\
	21	36.9947477599368\\
	22	53.8411746423851\\
	23	45.1937036301997\\
	24	30.9490730229206\\
	25	52.3758190306741\\
	26	248.903254078858\\
	27	81.0418432783008\\
	28	15.8533958622452\\
	29	106.470539168753\\
	30	15.523998402342\\
	31	37.8688514010372\\
	32	5.99857623213558\\
	33	5.8329380439343\\
	34	1.14538487796949\\
	35	0.954340898549118\\
	36	1.65266967614921\\
	37	4.18680295035514\\
	38	0.500363041598939\\
	39	0.409974848261882\\
	40	0.231888196382084\\
	41	7.13828835067339\\
	42	0.273536669043478\\
	43	1.1296097766348\\
	44	0.287836924798691\\
	45	0.968736006483894\\
	46	0.952382958141573\\
	47	0.159461609787263\\
	48	0.195157026360186\\
	49	0.425589337441317\\
	50	0.0934795157091343\\
	51	0.16210839078406\\
	52	0.054248611874899\\
	53	0.0440982574951139\\
	54	0.0245633223165308\\
	55	0.0663979261032342\\
	56	0.16378986381476\\
	57	0.0634482159182281\\
	58	0.0159324798318779\\
	59	0.00606656750867422\\
	60	0.000836121793702744\\
	61	2.37216196975627e-05\\
};
\addlegendentry{$\epsilon_{\text{est}}$ (proposed)}
\end{axis}
\end{tikzpicture}%
	\caption{Coax-fed dielectric resonator filter: results for Tests 2, 3 and 4: Convergence of Algorithm~\ref{alg:rom_errest}.}
	\label{fig:Brauer_Test234_conv}
\end{figure}
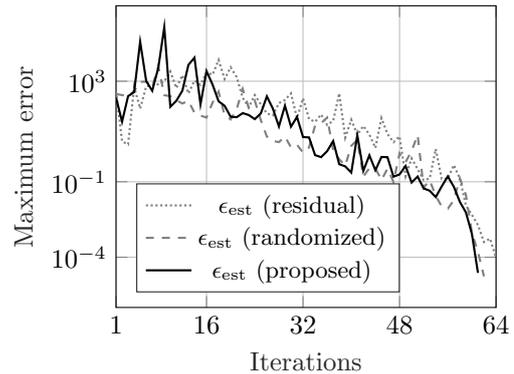

\setlength\fheight{3cm}
\setlength\fwidth{3cm}

\subsection{Example VI :  Inline Dielectric Resonator Filter}
\begin{figure}[t!]
	\centering
	\includegraphics[width=0.85\linewidth]{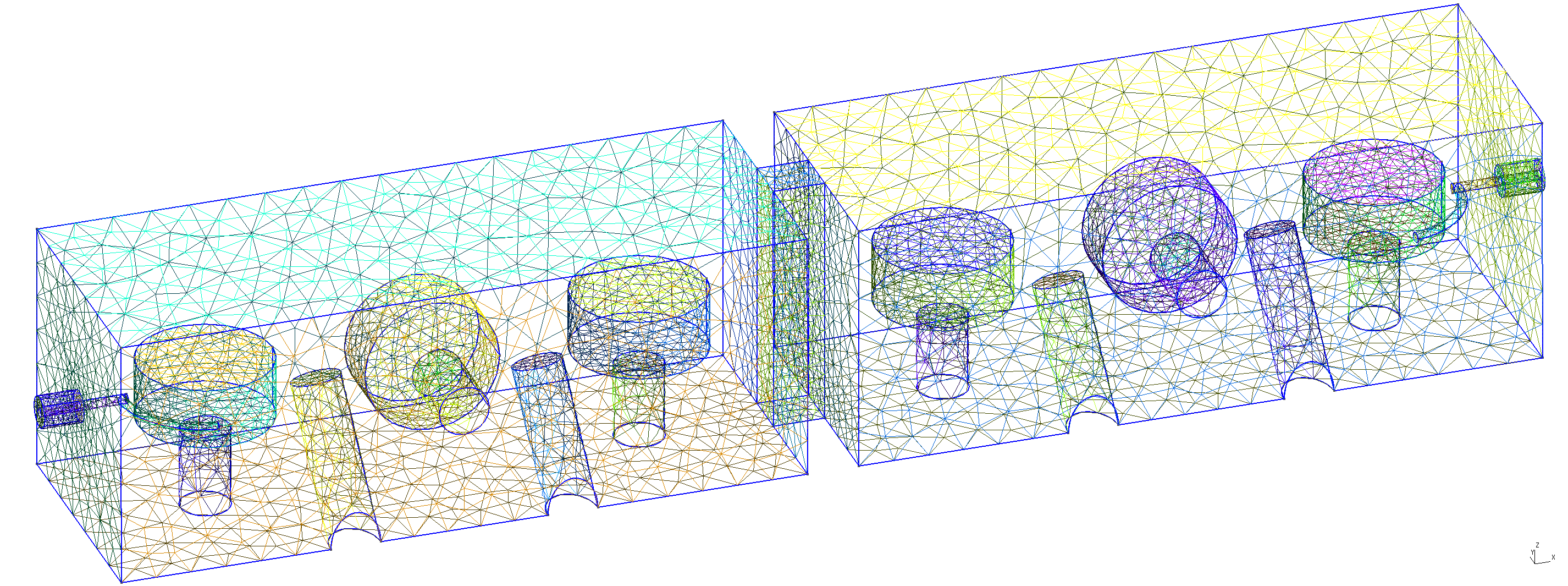} 
	\caption{Geometry of the sixth-order inline dielectric resonator filter from \cite{BasS12}.}
	\label{fig:Snyder_Geometry}
\end{figure}

The final example we discuss is the sixth-order inline dielectric resonator (I-DR) filter from \cite{BasS12} shown in Fig.~\ref{fig:Snyder_Geometry}. It is made of two cascaded triple-resonator configurations, resulting in six inline dielectric resonators. In addition to the frequency $f$, we consider the dielectric constants ($d_{1}, d_{2}$) from the two triple-resonators as additional parameters ($d_{1}$, $d_{2}$ stand for the relative dielectric permittivity of the first three and last three inline dielectric resonators, respectively), resulting in a three-parameter system. We have $f \in [2.10, 2.25]$ GHz, $d_{1}, d_{2} \in [76.5, 77.5]$. Further, the reference dielectric constant of the dielectric material is $d_{\text{ref}} = 77.0$. The discretized model has dimension $n=229,890$. The training set for this example is obtained using $1000$ samples of the frequency and $2$ samples each for the dielectric constants $d_{1}, d_{2}$. A cartesian grid of dimension $1000 \times 2 \times 2$ is formed resulting in a training set with $4000$ parameter samples. The ROM tolerance is $\texttt{tol} = 10^{-3}$.

Similar to the previous example, we do not show the true error convergence for this example in Tests 2 - 4 due to the huge computational costs.

%

\subsubsection{Test 2. Residual norm as a heuristic error estimator}
For the inline dielectric resonator filter, using the residual norm in Algorithm~\ref{alg:rom_errest} results in a ROM of dimension $r = 30$. As shown in Fig.~\ref{fig:IDRF_Tests234_conv} (dotted line), the greedy algorithm takes $15$ iterations to converge, taking around $3$ minutes and $50$ seconds.\\

%
\subsubsection{Test 3. Randomized error estimator}
We next apply the randomized error estimator $\tilde{\Delta}(\mu)$ from \cite{morSemZP18}. For this example, we sample $K=10$ random vectors each of dimension $n$ to obtain $V_{rd}$, with $\texttt{tol}_{rd} = 0.5$. We use the $\texttt{simdTwister}$ random number generator, with a seed of $10$. The total time to obtain $V_{rd}$ is roughly $206$ minutes. Fig.~\ref{fig:IDRF_Tests234_conv} (dashed line) illustrates the convergence of the greedy algorithm. It converges in $12$ iterations, taking $3$ minutes and $23$ seconds. The resulting ROM has dimension $r = 24$. Note that the dimension of this ROM is smaller than the one obtained from Test 2.
\\
%

\subsubsection{Test 4. Proposed error estimator}
Finally, we use the proposed error estimator $\| \tilde{e}(\mu) \|$ in Algorithm~\ref{alg:rom_errest}. The ROM obtained has dimension $r = 28$. Fig.~\ref{fig:IDRF_Tests234_conv} shows the convergence of the greedy algorithm, with the convergence achieved in $14$ iterations. When compared to the randomized error estimator, the proposed approach has a ROM of larger dimension. But, the time taken by the proposed error estimator is significantly shorter, at $8$ minutes and $5$ seconds. 

\begin{figure}[tbp]
	\centering
	\input{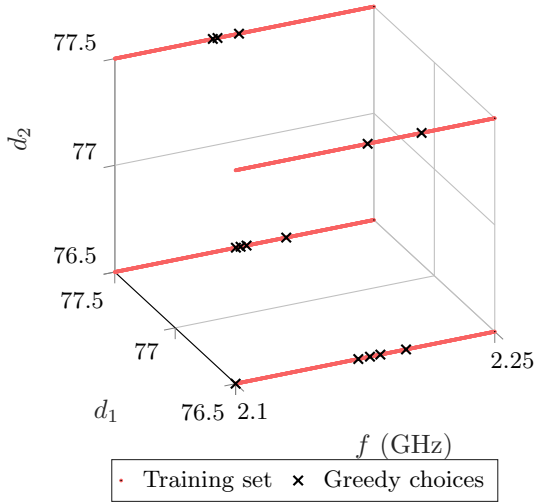}
	\caption{Inline dielectric resonator filter: greedy parameters $\mu^{*}$ selected by Algorithm~\ref{alg:rom_errest} using the proposed error estimator.}
	\label{fig:IDRF_Test4_gp}
\end{figure}

In Fig.~\ref{fig:IDRF_Test4_gp} we plot the locations of the greedy parameters $\mu^{*}$ for Test 4. It is worth noting that most of the samples are present in the frequency range $[2.15, 2.20]$ GHz where the resonances occur. Finally, we test the performance of the parametric ROM obtained from Test 4 on two parameter samples which were not considered in the training set $\Xi$. For this, we plot the magnitude and phase of the scattering parameter responses for two parameters $(d_{1}, d_{2}) = (76.6, 76.9)$ and $(d_{1}, d_{2}) = (77.2, 76.9)$. As seen from Fig.~\ref{fig:IDRF_Test4_sparam_test1} and Fig.~\ref{fig:IDRF_Test4_sparam_test2}, the scattering parameter response obtained from the ROM is nearly identical to the one obtained from the FOM solver. This highlights the robustness of the ROM generated using the \emph{inf-sup}-constant-free error estimator and also its ability to generalize well beyond the training data.\\
\setlength\fheight{4cm}
\setlength\fwidth{5cm}
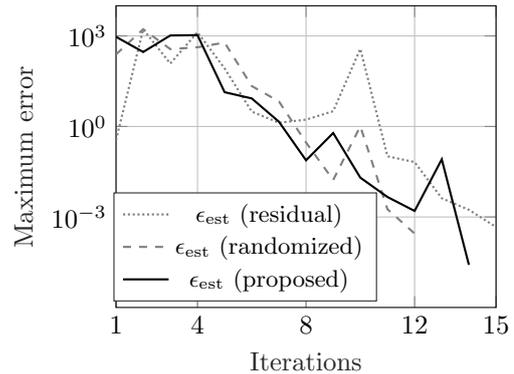
\begin{figure}[tbp]
	\centering
%
%
\definecolor{mycolor1}{rgb}{0.00000,0.44700,0.74100}%
\begin{tikzpicture}

\begin{axis}[%
width=\fwidth,
height=\fheight,
at={(0\fwidth,0\fheight)},
scale only axis,
xmin=1,
xmax=15,
xtick = {1, 4, 8, 12, 15},
xlabel style={font=\color{white!15!black}},
xlabel={Iterations},
ymode=log,
ymin=1e-06,
ymax=10000,
ytick = {1000,1,0.001},
yminorticks=true,
ylabel style={font=\color{white!15!black}},
ylabel={Maximum error},
axis background/.style={fill=white},
xmajorgrids,
ymajorgrids,
yminorgrids,
legend style={at={(0.34,0.02)}, anchor=south, legend cell align=center, align=left, font=\small, draw=white!15!black}
]
\addplot [color=gray, densely dotted, thick]
  table[row sep=crcr]{%
1	0.35697469899361\\
2	1657.38598700192\\
3	122.029895142884\\
4	1330.18681425104\\
5	84.1583647644299\\
6	3.11122518716958\\
7	1.32643738972739\\
8	1.68033953790163\\
9	3.20022750881931\\
10	375.138003762359\\
11	0.102622043829555\\
12	0.0659369116715242\\
13	0.00418019653610992\\
14	0.00172349397013389\\
15	0.000465483130717871\\
};
\addlegendentry{$\epsilon_{\text{est}}$ (residual)}
%
\addplot [color=gray,dashed,thick]
table[row sep=crcr]{%
	1	244.782566358488\\
	2	1673.22526848808\\
	3	366.121721105805\\
	4	422.125043705507\\
	5	602.185664539419\\
	6	21.8719385532036\\
	7	6.78143734401966\\
	8	0.293774439591825\\
	9	0.0168570405250034\\
	10	1.01424941895509\\
	11	0.00180611259313684\\
	12	0.00028208844689114\\
};
\addlegendentry{$\epsilon_{\text{est}}$ (randomized)}
%
\addplot [color=black, thick]
table[row sep=crcr]{%
	1	935.691403523452\\
	2	294.035906709638\\
	3	1039.21884034949\\
	4	1066.63291088971\\
	5	13.6889397849503\\
	6	8.4973699410753\\
	7	1.4380153253316\\
	8	0.075609282019707\\
	9	0.608248321695134\\
	10	0.0201087578650226\\
	11	0.00459427480560175\\
	12	0.00158361572503189\\
	13	0.0815030505630689\\
	14	2.59186538701597e-05\\
};
\addlegendentry{$\epsilon_{\text{est}}$ (proposed)}
\end{axis}
\end{tikzpicture}%
	\caption{Inline dielectric resonator filter: results for Tests 2, 3 and 4: Convergence of Algorithm~\ref{alg:rom_errest}.}
	\label{fig:IDRF_Tests234_conv}
\end{figure}


\setlength\fheight{5cm}
\setlength\fwidth{5cm}

\setlength\fheight{3cm}
\setlength\fwidth{3cm}

\begin{figure}[t!]
	\centering
	\input{figures/IDRF_sparam_FOMROM_MagPha_testparam1}
	\caption{Inline dielectric resonator filter: scattering parameter responses of FOM and ROM computed from Test 4 for parameter $(d_{1}, d_{2}) = (76.6, 76.9)$. Top: The magnitudes of the scattering parameter responses $|S_{11}|$, $|S_{21}|$; Bottom: The phases $\angle(S_{11}), \angle(S_{21})$.}
	\label{fig:IDRF_Test4_sparam_test1}
\end{figure}
\begin{figure}[t!]
	\centering
	\input{figures/IDRF_sparam_FOMROM_MagPha_testparam2}
	\caption{Inline dielectric resonator filter: scattering parameter responses of FOM and ROM computed from Test 4 for parameter $(d_{1}, d_{2}) = (77.2, 76.9)$. Top: The magnitudes of the scattering parameter responses $|S_{11}|$, $|S_{21}|$; Bottom: The phases $\angle(S_{11}), \angle(S_{21})$.}
	\label{fig:IDRF_Test4_sparam_test2}
\end{figure}


%

\begin{table}[t!]
	\centering
		\caption{Offline time to generate projection matrix $V$}
		\label{tab:offtime}
		\def\arraystretch{1.5}
		\begin{tabular}{|c|c|c|c|c|}
			\hline
			\multirow{2}{*}{Example}   &  \multicolumn{4}{c|}{Time Taken (mins)} \\ 
			\cline{2-5}		& \multicolumn{1}{c|}{Test 1} & \multicolumn{1}{c|}{Test 2} &  \multicolumn{1}{c|}{Test 3} & \multicolumn{1}{c|}{Test 4}\\
			\hline
			\hline
			Dual-mode filter           & $6.49$   & $0.21$    & $3.71$    &  $0.58$ \\ 
			\hline
			SIW antenna                &   -  &  $12.35$    &  $38.78$    &  $21.38$\\ 
			\hline
			AV antenna    			   &   -  &  $30.50$    &  $116.15$    &  $56.00$\\ 
			\hline
			DR antenna                 &   -  &  $48.88$   &  $159.25$    &  $90.23$ \\ 
			\hline
			C-DR filter                 &   -  &  $12.55$   &  $216.80$    &  $29.93$ \\
			\hline
			I-DR filter                 &   -  &  $3.82$   &  $209.31$    &  $8.10$ \\
			\hline
		\end{tabular}
		\def\arraystretch{1.5}
\end{table}











The time taken (in minutes) for the greedy algorithms for all the six examples are summarized in Table \ref{tab:offtime}. The greedy algorithm based on the residual norm requires the least computational time among the four error estimators. However, as illustrated in our numerical examples, the accuracy of the residual-norm based error estimator is not uniformly good. The residual norm overestimates the true error with an effectivity factor up to $50$. Moreover, although all the examples show that the residual norm gives ROMs with acceptable accuracy, this is not theoretically guaranteed \emph{a priori}, meaning it may not always give accurate ROMs. Furthermore, this situation may also lead to ROMs with larger-than-necessary dimension. The proposed error estimator, while acceptably more expensive than the residual estimator, gives a far superior estimation of the true error. In contrast to the residual norm, the proposed error estimator gives ROMs with theoretically guaranteed accuracy. This deserves the sacrifice of an extra but limited amount of offline cost, if a reliable ROM, a goal of MOR more important than the marginally increased offline time, is desired. Also, in our tests, it requires significantly shorter offline time when compared to the randomized error estimator.

%
\begin{table}[t!]
	\scriptsize
	\centering
	\caption{Comparison of computational costs: Brute-force vs. Proposed approach}
	\label{tab:cost_mor_vs_brute}
	\def\arraystretch{1.5}
	\begin{tabular}{|c|c|c|c|}
		\hline
		\multirow{2}{*}{Example}&\multicolumn{2}{c|}{Time Taken (s)}&\multirow{2}{*}{Speedup}\\ 
		\cline{2-3}		& \multicolumn{1}{c|}{Brute-force} & \multicolumn{1}{c|}{\small{offline} (\scriptsize{Alg.~\ref{alg:rom_errest}}) \small{+ online}} & \\
		\hline
		\hline
		Dual-mode filter           & 143    & 35 + 0.008    & 4.1\\ 
		\hline
		SIW antenna                & 9273  &  1283 + 0.007    &7.2\\ 
		\hline
		AV antenna    			   & 9110  &  3360 +  0.018   & 2.7\\ 
		\hline
		DR antenna                 & 24843  &  5414 + 0.015    & 4.6\\ 
		\hline
		C-DR filter                 & 4600  &  1796 + 0.05    & 2.6\\ 
		\hline
		I-DR filter                 & 6271  &  486 + 0.015    & 12.9\\ 
		\hline
	\end{tabular}
	\def\arraystretch{1.5}
\end{table}
\begin{table}[t!]
	\scriptsize
	\centering
	\caption{Maximum error over the test parameter set for Tests 2, 3 and 4}
	\label{tab:test_errors}
	\def\arraystretch{1.5}
	\begin{tabular}{|c|c|c|c|c|}
		\hline
		\multirow{2}{*}{Example}&\multirow{2}{*}{$\texttt{tol}$}&\multicolumn{3}{c|}{Maximum error over test set}\\ 
		\cline{3-5}		& &\multicolumn{1}{c|}{Test 2} & \multicolumn{1}{c|}{Test 3} & \multicolumn{1}{c|}{Test 4}\\
		\hline
		\hline
		Dual-mode filter           & $10^{-6}$ &$5.22\cdot 10^{-8}$    & $4.19\cdot 10^{-7}$    & $6.32\cdot 10^{-7}$\\ 
		\hline
		SIW antenna                & $10^{-4}$ & $2.26 \cdot 10^{-6}$  &  $3.33 \cdot 10^{-5}$    &$2.20 \cdot 10^{-5}$\\ 
		\hline
		AV antenna    			   & $10^{-3}$ & $1.06 \cdot 10^{-4}$  &  $8.00 \cdot 10^{-4}$   & $8.36 \cdot 10^{-4}$\\ 
		\hline
		DR antenna                 &$10^{-3}$ & $1.80 \cdot 10^{-5}$  &  $4.04 \cdot 10^{-4}$    & $2.13 \cdot 10^{-4}$\\ 
		\hline
		C-DR filter                 &$10^{-4}$ &  $6.98 \cdot 10^{-8}$  &  $5.70\cdot 10^{-5}$  & $5.80 \cdot 10^{-5}$\\ 
		\hline
		I-DR filter                 &$10^{-3}$ & $4.96 \cdot 10^{-5}$  &  $5.52 \cdot 10^{-5}$    & $5.05 \cdot 10^{-5}$\\ 
		\hline
	\end{tabular}
	\def\arraystretch{1.5}
\end{table}

In Table~\ref{tab:cost_mor_vs_brute}, we compare the total computational cost of a reduced-order model generated using the proposed approach (Test 4) with that of a brute-force FOM simulation, over a test set of parameters. For all the six examples, the test set consists of parameter samples not used in the offline stage of MOR. For the filter model, the test set consists of 200 samples, while for the three antenna examples, the test set has 100 samples in their corresponding parameter domains. For the last two examples, viz., the C-DR and I-DR filters, due to their larger storage requirements, we consider a smaller test set in comparison to the training set. In case of the C-DR filter, the test set contains $500$ new frequency locations. Since the I-DR filter is parametric, we consider $4$ new parameter combinations $(d_{1}, d_{2})$ not present in the training set. Each combination corresponds to $100$ frequency samples. The offline time reported for Algorithm~\ref{alg:rom_errest} is the same as in Table~\ref{tab:offtime}, but in seconds. It is evident that our proposed method has significantly reduced the computational time, especially for models with large size. The comparatively lower speedup for the dual-mode filter example is owing to its relatively smaller dimension $n = 36,426$ such that the solver requires less time to compute the solution. Moreover, the offline time of the proposed method is independent of the number of testing samples, and its online time is almost negligible compared with the offline time. Consequently, when the amount of testing samples is increased, the speedup over the brute-force approach will be even higher. 

In Table~\ref{tab:test_errors}, we show the maximum error over the test set, for each of the six examples corresponding to Tests 2, 3, and 4. We notice that that the maximum test errors for all the examples satisfy the defined tolerance. It can also be seen that the lowest errors result from Test 2. This is mainly due to the fact that the ROM obtained from Test 2 has the largest dimension. This also indicates that using the residual as the error estimator leads to overestimating the true error for those examples. The maximum error over the test set resulting from Tests 3 and 4 are very close to each other and to the assigned tolerance. This again justifies the effectivity of the corresponding error estimators. The error estimators can accurately estimate the true error with effectivity ($\texttt{eff}$ in the figures) close to $1$, thus producing ROMs with true errors also close to but below the tolerance. However, the computational cost of Test 3 is significantly larger. Thus, we can see that the proposed error estimator used in Test 4 offers a good compromise between the dimension of the ROM and its ability to offer good performance for unseen parameter values.

Finally, for frequency-domain simulation of problems with, e.g., geometrical or material parameters, which queries not only many frequency samples but also many parameter samples, further  efficiency will be achieved by using the proposed method as can be seen from the significant speedup achieved for the I-DR filter which is a three-parameter example.
\section{Conclusions}
\label{sec:conclusions}
We have introduced a novel \textit{a posteriori} error estimator capable of accurately estimating the state error, even for systems where the \emph{inf-sup} constant is very small. Many engineering systems involving resonant electromagnetic devices, such as microwave filters, antennas show this behaviour. A waveguide filter model, models of three types of antenna and two different types of dielectric resonator filters are used to demonstrate the robustness of the proposed error estimator. The proposed error estimator outperforms the standard estimator and a recently proposed one in the literature, both theoretically and numerically, thus showing its great potential in electromagnetic simulation and analysis. Although the proposed error estimator needs more offline time than the residual-based approach, it derives a ROM with guaranteed accuracy and is, therefore, more reliable. As a result, compact reduced order models for challenging real-life applications have been obtained.


\section*{Acknowledgments}%
\addcontentsline{toc}{section}{Acknowledgments}

S. Chellappa was supported by the International Max Planck Research School for Advanced Methods in Process and Systems Engineering (IMPRS-ProEng) when this work was done.


\addcontentsline{toc}{section}{References}
\bibliographystyle{plainurl}
\bibliography{InfSup_Free_Error_Estimation}
  
\end{document}